\documentclass[12pt]{amsart}
\usepackage{amsthm}
\usepackage{amssymb}
\usepackage{amsmath}
\usepackage{xcolor}

\usepackage{indentfirst}
\theoremstyle{definition}
\newtheorem{prop}{Proposition}[section]
\newtheorem{theo}[prop]{Theorem}
\newtheorem{defi}[prop]{Definition}
\newtheorem{ex}[prop]{Example}
\newtheorem{coro}[prop]{Corollary}
\newtheorem{lemma}[prop]{Lemma}
\newtheorem{rem}[prop]{\textbf{Remark}}
\usepackage[all]{xy}
\usepackage{enumerate}
\usepackage[a4paper,centering,left=2cm,right=2cm,top=2.3cm,bottom=2.3cm]{geometry}
\usepackage{graphicx}
\usepackage{float}
\usepackage{bm}
\usepackage{color}
\usepackage{mathtools}
\usepackage[style=alphabetic]{biblatex}
\usepackage{hyperref}
\usepackage{microtype}
\usepackage{colonequals}
\usepackage{verbatim}
\usepackage{microtype}
\usepackage{tikz}
\usepackage{mathtools}
\usepackage{comment}
\usepackage[normalem]{ulem}
\title{N\'eron Functions on Jacobians of
Odd-Degree Hyperelliptic Curves}
\author{Tianci Kang}
\address{T.Kang, 
Bernoulli Institute,
Rijksuniversiteit Groningen,
Nijenborgh 9,
9747 AG Groningen,
The Netherlands}
\email{t.kang@rug.nl}
\date{}
\begin{document}
\begin{abstract}
We express N\'eron functions on Jacobians of odd-degree hyperelliptic curves of arbitrary genus in terms of admissible pairings, using the pullback of the theta divisor via a generalized Abel--Jacobi morphism. Over non-archimedean local fields, we relate N\'eron functions to local N\'eron pairings on the curve and, moreover, deduce a comparison between unramified $p$-adic N\'eron functions and unramified local Coleman--Gross height functions. Furthermore, we decompose the N\'eron function into an intersection term $I$ and a Green function term $T$, and compute $T$ for specific semistable reduction types. We compare our decomposition with other decompositions in the literature.     
\end{abstract}
\maketitle
\tableofcontents
\section{Introduction}
N\'eron \cite{Neron1965quasi} and Tate constructed the canonical height function on an abelian variety $A$ over a number field $F$ associated to a symmetric line bundle $L$. It can be expressed as a sum of local components, called \textit{N\'eron functions}, defined over the completions of $F$. Algorithms for computing canonical heights have numerous applications, for instance, in computing generators of $A(F)$ and the regulator of $A/F$, which appears in the Birch--Swinnerton--Dyer conjecture. 

Explicit formulas for N\'eron functions on elliptic curves have been described in \cite[Chapter~VI]{SilvermanAdvancedTopics}. For an abelian variety, one of the ways to construct N\'eron functions is using an explicit embedding of the Kummer variety associated to the abelian variety. Explicit embeddings have been described for Jacobians of genus $2$ curves, see \cite{Flynncanonical}, and for Jacobians of genus $3$ hyperelliptic curves, see \cite{Stollgenus3}. Using Kummer embeddings, algorithms for computing N\'eron functions on Jacobians of hyperelliptic curves of genus $g\leq3$ are obtained. For hyperelliptic curves of genus $4$, this is work in progress by Ludwig Fürst. However, the dimension of the embedding projective space is $2^g-1$, which makes such explicit embeddings difficult to control as $g$ increases. 

Let $J(C)$ be the Jacobian of a curve $C\colon y^2=f(x)$ of genus $2$ over a non-archimedean local field $K$ of characteristic different from $2$, where $f$ is a polynomial of degree $5$. Let $\Theta$ be the theta divisor associated to the point at infinity $\infty$. The N\'eron function $\lambda_\Theta$ on $J(C)(K)\setminus|\Theta|$ can be expressed in terms of Zhang's admissible pairings on $C$ via the difference morphism $C\times C\to J$, see \cite[Proposition~6.5]{Bianchicoleman}. In this paper, we obtain an explicit formula for the N\'eron function $\lambda_\Theta$ in terms of admissible pairings on an odd-degree hyperelliptic curve $C$ of arbitrary genus $g$ over local fields of characteristic different from $2$, see Theorem~\ref{Thm6.2}. The formula in Theorem~\ref{Thm6.2} does not require any explicit projective embedding of the Jacobian or the Kummer variety. In particular, over non-archimedean local fields, we show that \cite[Proposition~6.5]{Bianchicoleman} is a special case of Theorem~\ref{Thm6.2}. The key step in the proof of Theorem~\ref{Thm6.2} is the pullback of the theta line bundle via a generalized Abel--Jacobi morphism, which is proved in Theorem~\ref{theo1n}. When $g=2$, the formula in Theorem~\ref{theo1n} recovers \cite[Lemma~4.10.2]{Moretmetriquespermises}, which determines the pullback of the theta line bundle via the difference map.  

For a hyperelliptic curve $C$ over a local field $K$, another local height function $\lambda_N\colon J(C)(K)\setminus|\Theta|\to\mathbb{R}$ is defined in \S\ref{ch7} using an extension of the N\'eron pairing on $C$ due to \cite{Grosslocalheights}. When $K$ is a non-archimedean local field, it is proved in Theorem~\ref{main} that the N\'eron function $\lambda_{2\Theta}$ is equal to the function $\lambda_N$ up to a simple correction term. This result has the following application to the theory of $p$-adic heights.

Let $A$ be an abelian variety over a number field $F$, and let $D\in\operatorname{Div}(A)$. For a prime number $p$ and a finite place $v$ of $F$, the $p$-adic N\'eron function $\lambda_{D,\chi_v}\colon A(F_v)\setminus|D|\to\mathbb{Q}_p$ with respect to $D$ and a local homomorphism $\chi_v\colon F_v^{\times}\to\mathbb{Q}_p$ is defined in \cite[\S4]{Bianchicoleman}, and is a local component of the canonical Mazur--Tate $p$-adic height, see \cite{MazurTateCanonicalheight}. When $\chi_v$ is unramified, the N\'eron function $\lambda_D$ and the $p$-adic N\'eron function $\lambda_{D,\chi_v}$ are scalar multiples of the same $\mathbb{Q}$-valued function. Let $C$ be a curve over $F$. The local Coleman--Gross height pairing $h^{\operatorname{CG}}_{\chi_v}\colon \operatorname{Div}^{0}(C_{F_v})\times\operatorname{Div}^{0}(C _{F_v})\to\mathbb{Q}_p$ with respect to $\chi_v$ is defined in \cite{ColemanGrossp-adic} for disjoint divisors and extended in \cite{BalakrishnanBesser} to arbitrary divisors. When $\chi_v$ is unramified, both the extended local Coleman--Gross height pairing and the extended N\'eron pairing are the same $\mathbb{Q}$-valued function up to a constant factor, see \cite[Proposition~3.4]{BalakrishnanBesser}. The local Coleman--Gross height function $\lambda^{\operatorname{CG}}_{\chi_v}\colon J(C)(K)\setminus|\Theta|\to\mathbb{Q}_p$ is defined by $\lambda^{\operatorname{CG}}_{\chi_v}(a)\coloneq h^{\operatorname{CG}}_{\chi_v}(D,D)$ where $D\in\operatorname{Div}^0(C _{F_v})$ is a certain representative of $a$, see \cite[Definition~6.1]{Bianchicoleman}.

When $v$ is a finite place of $F$ and $C/F_v$ is a curve of genus $2$, it is proved in \cite[Theorem~6.3]{Bianchicoleman} and \cite[Theorem~6.15]{Bianchicoleman} that the $p$-adic N\'eron function $\lambda_{2\Theta,\chi_v}$ on $J(C)(K)\setminus|\Theta|$ is equal to the local Coleman--Gross height function $\lambda^{\operatorname{CG}}_{\chi_v}$ up to a constant factor. When $C/F_v$ is a hyperelliptic curve of genus $g$, Theorem~\ref{main} generalizes \cite[Theorem~6.3]{Bianchicoleman} by showing that the unramified $p$-adic N\'eron function $\lambda_{2\Theta,\chi_v}$ differs from the unramified local Coleman--Gross height function $\lambda^{\operatorname{CG}}_{\chi_v}$ by a simple correction term.

Furthermore, Theorem~\ref{main} can be interpreted as the local counterpart of the formula $\sum_v\langle D,E\rangle_{v,N}=-h_{2\Theta}([D],[E])$ due to Faltings and Hriljac, which relates the sum of N\'eron pairings on a curve to the canonical height on the Jacobian of a curve. It is natural to conjecture that Theorem~\ref{main} also holds in the archimedean setting. The absence of explicit formulas for the Arakelov--Green function is the main obstacle to proving the conjecture. Remark~\ref{conjecture} shows that the conjecture holds when $F$ is a number field with only one archimedean place.    

It is shown by \cite[Theorem~4.6]{Zhangadmissiblepairing} that when $C$ is a curve with semistable reduction over a non-archimedean local field, the admissible pairing can be expressed using intersection theory on a regular model of $C$ and harmonic analysis on the reduction graph of $C$. This description can be used to decompose the N\'eron function $\lambda_\Theta$. 

We study further properties of our construction of N\'eron functions and compare it to other constructions in the literature. In Corollary~\ref{decomposition}, the N\'eron function $\lambda_\Theta$ is decomposed into an intersection multiplicity function $I$ and a Green function term $T$. 

Let $A$ be an abelian variety over a non-archimedean local field $K$, and let $D$ be a divisor on $A$. It is shown in \cite[\S III.4, Th\'eor\`eme 1, p. 319]{Neron1965quasi} and \cite[Theorem~5.1]{LangDiophantine} that there exists a map $\gamma_D\colon\Phi\to\mathbb{R}$ such that the N\'eron function associated to $D$ can be written as \[\lambda_D(a)=i(a,D)+\gamma_D(\textbf{a})\] where $a\in A(K)\setminus|D|$ and $\textbf{a}$ is the specialization of $a$ in the component group $\Phi$ of the N\'eron model of $A$. However, it is quite difficult to obtain an explicit formula for $\gamma_D$. When $D$ is symmetric and ample and defines a principal polarization, and $A$ has semi-abelian reduction, it is shown in \cite[Theorem~B]{dejongLocalheightsberkovich} that $\gamma_D$ can be expressed in terms of a certain tropical theta characteristic of $D$.

Using the Kummer surface, a correction term $\mu$ for N\'eron functions on $J(C)$ is defined in \cite[Definition~3.1]{MullerCanonical} when $C$ is a curve of genus $2$ over a non-archimedean local field $K$. It is proved in \cite[Theorem~7.4]{MullerCanonical} that $\mu$ factors through $\Phi$ whenever the given integral Weierstrass model of $C$ has rational singularities. When $\mu$ factors through $\Phi$, we can prove that $\mu$ is equal to $-2T$ up to a constant. Moreover, \cite[\S4.4]{DonkerCanonical} computes the correction term $\gamma_\Theta$ for curves of genus $2$ with various semistable reduction types by choosing a specific tropical theta characteristic in \cite[Theorem~B]{dejongLocalheightsberkovich} and assuming that the chosen theta characteristic gives rise to $\gamma_\Theta$. Furthermore, with the same choice of the tropical theta
characteristic, $\mu$ and $\gamma_\Theta$ are related in \cite[Proposition~4.5.13]{DonkerCanonical} for certain reduction types. Assuming the choice of the tropical theta
characteristic is correct, it was shown that $\mu=-2\gamma_\Theta$ when $\mu$ factors through $\Phi$.  

Moreover, given the semistable reduction type of a hyperelliptic curve $C$ of genus $g$, an explicit formula for $T$ can be obtained from the formula for Zhang's Arakelov--Green function with respect to the canonical divisor in Lemma~\ref{kayagreen}. Hence, we can compute all possible values of $T$ and determine whether $T$ factors through $\Phi$ by considering all possible specializations of points on $C$. In \S\ref{ch8.2}, we carry out this computation in the following two cases of hyperelliptic curves with semistable reduction: a hyperelliptic curve of genus $3$ such that the special fiber of the stable model is an irreducible curve with a node of thickness $3$, whose normalization has genus $2$; and a curve of genus $2$ such that the special fiber of the stable model consists of two components of genus $1$ meeting at a node of thickness $2$. We prove that $T$ does not factor through $\Phi$ in either case. 

The paper is organized as follows. In \S\ref{ch2}, we introduce N\'eron pairings on curves. In \S\ref{ch3}, we introduce admissible pairings on curves and Zhang's construction of admissible pairings in terms of intersection multiplicities and Arakelov--Green functions on the reduction graphs of curves over non-archimedean local fields. In \S\ref{ch5}, we define the generalized Abel--Jacobi morphism and the difference morphism on curves, and compute the pullback of the theta line bundle via these two different morphisms. In \S\ref{ch6}, an explicit formula for N\'eron functions on Jacobians of odd-degree hyperelliptic curves over a local field is given in Theorem~\ref{Thm6.2} in terms of admissible pairings. In \S\ref{ch7}, Theorem~\ref{main} establishes a relation of N\'eron functions with N\'eron pairings on odd-degree hyperelliptic curves of arbitrary genus. Its $p$-adic analogue generalizes \cite[Theorem~6.3]{Bianchicoleman}. In \S\ref{ch8}, we introduce an explicit formula for the Green function given by \cite[Equation 4]{KayaComputation} and use it to compute $T$ in some examples.                            

\textbf{Notation.}  

Throughout this paper, a curve of genus $g$ over a field $K$ is always projective, smooth, and geometrically irreducible. For a divisor $D$ of degree zero on a curve $C/K$, we denote its class in $J(C)=\operatorname{Pic}^0(C)$ by $[D]$. 

Let $\alpha\in C(K)$ be a rational point. The theta divisor $\Theta^\alpha$ with respect to $\alpha$ on $J(C)$ is the prime divisor whose support is \[\{[D-(g-1)\alpha]\,|\,D\ge0,\deg D=g-1\}.\] The line bundle $\mathcal{O}_J(\Theta^\alpha)$ associated to the divisor $\Theta^\alpha$ is denoted by $\theta^\alpha$.        

An odd-degree hyperelliptic curve $C$ of genus $g\ge1$ over a field $K$ of characteristic different from $2$ is a smooth projective curve with affine model\[y^2=f(x)\]
where $f\in K[x]$ is a separable polynomial of degree $2g+1$. There is a unique point $\infty\in C(K)$ on $C$ at infinity. The hyperelliptic involution on $C$ is denoted by $\iota$. Throughout the paper, all hyperelliptic curves are assumed to be of odd degree over a field $K$ of characteristic different from $2$. 

For functions $A$ and $B$ with values in an abelian group, the notation $A\equiv B$ means that $A$ is equal to $B$ up to an additive constant.

\textbf{Statement on AI usage.}  

During the writing of this paper, the author only asked ChatGPT to check the manuscript for grammatical and mathematical errors. The author did not use AI tools to develop any original mathematical arguments.

\textbf{Acknowledgements.}
I would like to express my sincere gratitude to my supervisor, Steffen M\"uller, for introducing me to this interesting problem. I have benefited a lot from our discussions. This paper would not have been possible without his insightful and detailed comments throughout the writing process. I would also like to thank Robin de Jong and Klaus K\"unnemann for their valuable comments on the first draft of this paper. I am grateful to Martin Lüdtke for helping me prepare a presentation based on the results of this paper at the Intercity Seminar on Arakelov Geometry 2026. I would like to thank Farbod Shokrieh, Marc Hindry and Paul Alexander Helminck for their thoughtful comments on my presentation at the seminar. Moreover, I would like to thank Quinten Donker for sharing his master's thesis and for our discussion on possible directions for further research.    

\section{N\'eron pairings on curves and N\'eron functions on abelian varieties}\label{ch2}

In this section, we briefly recall N\'eron pairings on curves. We define the N\'eron pairing for divisors with disjoint support and then discuss its extension to divisors with common support.

\subsection{N\'eron pairings for divisors with disjoint support on curves}
\begin{theo}\cite[Theorem 1.1]{Hriljacheights}\label{thm2.1}
Let $C$ be a curve over a field $K$ with a nontrivial absolute value $v$. Assume that $C(K)$ is Zariski dense in $C$. There is a unique $\mathbb{R}$-valued pairing $\langle D,E\rangle_N^K$ on divisors $D$ and $E\in\operatorname{Div}^0(C(K))$ with disjoint support, satisfying the following properties.
    \begin{enumerate}
        \item The pairing is symmetric and bilinear.
        \item If $D=\operatorname{div}(f)$ is principal, then $\langle D,E\rangle_N^K=-\log|f(E)|_v$.
        \item Let $x_0\in C(K)\setminus|D|$, then the map $C(K)\setminus|D|\to\mathbb{R}$ defined by $x\mapsto\langle D,(x-x_0)\rangle_N^K$ is continuous. 
    \end{enumerate}
\end{theo}
\begin{rem}\cite[\S2]{Grosslocalheights}\label{extendedGross}
There is an extension of the pairing $\langle D,E\rangle_N^K$ such that it is independent of the base field. Let $H$ be a finite extension of $K$, and let $D$ and $E$ be relatively prime divisors of degree zero over $K$. It is easy to verify that \[\langle D,E\rangle^H_N=[H:K]\langle D,E\rangle_N^K.\] Then for any finite field extension $H/K$ such that $D$ and $E$ are pointwise rational over $H$, we define the general N\'eron pairing as follows\[\langle D,E\rangle_N\coloneq\frac{1}{[H:K]}\langle D,E\rangle_N^H.\] Throughout the paper, $\langle D,E\rangle_N$ always denotes the absolute N\'eron pairing.       
\end{rem}

When $C$ is a curve over a non-archimedean local field $K$, the N\'eron pairing on $C$ can be interpreted using intersection multiplicities on arithmetic surfaces, see \cite[\S3]{Grosslocalheights} and \cite[III.5 Theorem~5.2]{LangArakelov}, as we now recall.       

Let $\mathcal{C}/\operatorname{Spec}\mathcal{O}_K$ be a flat proper regular model for $C$. The group of rational fibral divisors on $\mathcal{C}$ is denoted by $\mathbb{Q}\operatorname{Div}_s(\mathcal{C})$. Let $\mathbb{Q}\mathcal{C}_s$ denote the $\mathbb{Q}$-vector space generated by the entire special fiber $\mathcal{C}_s$. We relate divisors of degree zero on $C$ to rational fibral divisors on $\mathcal{C}$ as follows.
\begin{prop}\cite[Theorem~1.3]{Hriljacheights}\cite[III Theorem~3.6]{LangArakelov}\label{othogonal}
    There exists a unique linear map \[\Phi\colon\operatorname{Div}^0(C(K))\to\mathbb{Q}\operatorname{Div}_s(\mathcal{C})/\mathbb{Q}\mathcal{C}_s\] such that for any $D\in\operatorname{Div}^0(C(K))$, the divisor $\bar{D}+\Phi(D)$ is orthogonal to $\operatorname{Div}_s(\mathcal{C})$.
\end{prop}
The following theorem due to Faltings and Hriljac shows that the N\'eron pairing on $C$ can be expressed in terms of the intersection multiplicity on $\mathcal{C}$.
\begin{theo}\cite[Theorem~1.6]{Hriljacheights}\label{intersection}
     For any divisors $D$ and $E$ in $\operatorname{Div}^0(C)$ with disjoint support, we have \[\langle D,E\rangle_N=\log q_v\cdot(\bar{D}+\Phi(D)\cdot\bar{E}+\Phi(E)),\] where $q_v$ is the norm of $v$.
\end{theo}
\subsection{N\'eron pairings for divisors with common support on curves}\label{ch2.2}
Gross extended N\'eron pairings to divisors with common support, see \cite[\S5]{Grosslocalheights}. Let $C$ be a curve over a non-archimedean local field $K$ with discrete valuation $v$. For any $K$-rational point $x\in C(K)$, fix a basis vector $\partial/\partial t_x$ of the tangent space at $x$ and choose a normalized parameter $z_x$ at $x$ such that $\partial z_x/\partial t_x=1$. A rational function $f\in K(C)^\times$ has a well-defined value at $x$ with respect to $t_x$ as follows\[f(x,t_x)\coloneqq\frac{f}{z_x^{\operatorname{ord}_xf}}(x).\]
It is easy to verify that $f(x,t_x)$ depends only on $t_x$ and not on $z_x$.
\begin{defi}\label{extended}
    Let $D,E\in\operatorname{Div}^0(C(K))$, and let $f\in K(C)^\times$ such that $D+\operatorname{div}(f)$ and $E$ are divisors with disjoint support. Fix a set-theoretical $K$-rational section $t$ of the tangent bundle on $C$, which assigns a nonzero tangent vector $t_P$ to every rational point $P$. The value of $f$ on the divisor $E=\sum_Pn_P(P)$ with respect to $t$ is defined as $f(D,t)\coloneqq\prod_Pf(P,t_P)^{n_P}$. 
    
    Then, the \textit{extended N\'eron pairing} between $D,E$ with respect to $t$ is defined as  \[\langle D,E\rangle_{N,t}\coloneqq\langle D+\operatorname{div}f,E\rangle_N+\log|f(E,t)|_v.\] Furthermore, the extended N\'eron pairing on $\operatorname{Div}^0(C)$ is defined by base change.  
\end{defi}
\begin{theo}\cite[\S5]{Grosslocalheights}
The extended N\'eron pairing satisfies the following properties:
\begin{enumerate}
        \item The pairing $\langle D,E\rangle_{N,t}$ does not depend on $f$.
        \item The pairing $\langle D,E\rangle_{N,t}$ depends on the choice of $t$, but it is unchanged if $\partial/\partial t$ is multiplied by $c\in K^\times$ with $|c|_v=1$.
        \item We have the identity $\langle \operatorname{div}f,E\rangle_{N,t}=-\log|f(E,t)|_v$ for $f\in K(C)^\times$.
\end{enumerate}
\end{theo}
Similarly to Theorem~\ref{intersection}, the extended N\'eron pairing has an interpretation in terms of intersection theory on arithmetic surfaces. First, we define the self-intersection of $K$-rational points. 
\begin{defi}
Let $t$ be a set-theoretical $K$-rational section of the tangent bundle of $C$. Let $P\in C(K)$, and let $\alpha\in K^\times$ such that $\alpha t(P)$ generates the tangent space of $\mathcal{C}$ at $P$. Let $\bar{P}$ be the Zariski closure of $P$ in $\mathcal{C}$. We define the \textit{self-intersection} of $\bar{P}$ with respect to $t$ as follows  
\[(\bar{P}\cdot\bar{P})_{t}\coloneqq-\operatorname{ord}_v(\alpha).\]
\end{defi} 
It is natural to extend $(\,\cdot\,)_t$ linearly to divisors on $\mathcal{C}$ whose restrictions to the generic fiber $C$ are pointwise $K$-rational. Furthermore, we have the following adjunction formula.
\begin{lemma}\label{adjunction}\cite[Lemma~3.11]{Bianchicoleman}
    Let $\omega$ be a meromorphic Kähler differential form on $C$. Let $P\in C(K)$ such that $P\notin|\operatorname{div}(\omega)|$. Let $t$ be the rational section determined by $\omega$ by duality. We have \[(\bar{P}\cdot\bar{P})_{t}+(\bar{P}\cdot\operatorname{div}(\bar{\omega}))=0\] where $\bar{\omega}$ is the differential form on $\mathcal{C}$ induced by $\omega$. 
\end{lemma}
More generally, the pairing extends to arbitrary divisors in $\operatorname{Div}^0(C)$ by base change. This generalization and the following result are discussed in \cite[\S3.2.1]{Bianchicoleman} for the unramified local Coleman--Gross height pairing, which agrees with the N\'eron pairing up to a multiplicative constant, see \S\ref{ch2.4}.  
\begin{prop}    
\cite[\S2, \S5]{Grosslocalheights}\label{extendedpairing}
Fix a set-theoretical rational section $t$ of the tangent bundle on $C$. For any two divisors $D,E\in\operatorname{Div}^0(C)$, we have \[\langle D,E\rangle_{N,t}=\log q_v\cdot\alpha_t(D,E),\] where $q_v$ is the norm of $v$ and $\alpha_t(D,E)\coloneq(\bar{D}+\Phi(D)\cdot\bar{E}+\Phi(E))\in\mathbb{Q}$. 
\end{prop}

\section{Admissible pairings over non-archimedean fields}\label{ch3}
Let $C$ be a curve over a local field $K$. Zhang constructed the admissible pairing of divisors on $\bar{C}=C_{\bar{K}}$ with disjoint support, see \cite[Theorem~4.6]{Zhangadmissiblepairing}. When $C$ is a curve over $\mathbb{R}$ or $\mathbb{C}$, the pairing is constructed using Arakelov theory, see \cite[\S4.5]{Zhangadmissiblepairing}. When $C$ is a curve over a non-archimedean local field, Zhang explicitly constructed the pairing using harmonic analysis on metrized graphs. Alternatively, the existence of the pairing follows from the bi-admissible metric on $\mathcal{O}_{C^2}(\Delta)$, obtained using admissible metrics on $J(C)$ and on $C$. The details are presented in this section. Yuan unified both archimedean and non-archimedean cases in terms of integrable metrics on $\omega_{C/K}$ and $\mathcal{O}_{C^2}(\Delta)$ in \cite[Theorem~A.1]{Yuanbigness}.

\begin{theo}
\cite[Theorem~4.6]{Zhangadmissiblepairing}\label{admissible}
Let $K$ be a local field, and let $\bar{K}$ be an algebraic closure of $K$. Let $C$ be a curve of genus $g\ge1$ defined over $K$, and let $\bar{C}=C_{\bar{K}}$ be the base change of $C$. Then there is the admissible pairing $\langle\;,\;\rangle_{\operatorname{ad}}$ defined on divisors in $\operatorname{Div}(\bar{C})$ with disjoint support, which satisfies the following properties:
\begin{enumerate}
\item The pairing $\langle\;,\;\rangle_{\operatorname{ad}}$ is symmetric and bilinear. 
\item For an open subset $U$ of $\bar{C}$ and a divisor $D\in\operatorname{Div}(\bar{C})$, let $f$ be a rational function such that $\operatorname{div}(f)|_U=D|_U$. Then $\langle D,x\rangle_{\operatorname{ad}}+\log|f(x)|$ is a bounded continuous function on $U$. 
\item For any nonzero rational function $g$ on $\bar{C}$, there exists a constant $c_g$ such that \[\langle\operatorname{div}(g),x\rangle_{\operatorname{ad}}=-\log|g(x)|+c_g\] for any $x\in\bar{C}\setminus|\operatorname{div}(g)|$.
\item Let $\alpha$ be a nonzero rational $1$-form on $\bar{C}$, and let $\mathcal{K}$ be the canonical divisor on $\bar{C}$. There exists a constant $c_\alpha$ such that for any $x\notin |\operatorname{div}(\alpha)|$ and any rational function $f$ over $\bar{C}$ with a simple pole at $x$, the following equation holds\[c_\alpha+\lim_{y\to x}\left(\langle \mathcal{K}+x,y\rangle_{\operatorname{ad}}+\log|f\alpha/df|(y)\right)=0.\] 
\end{enumerate}   
\end{theo}

\subsection{Admissible metrics }
Let $X$ be a projective variety over a local field $K$. Let $(L,||\cdot||)$ be a continuous and bounded metrized line bundle on $X$ (see \cite[Definition~2.7.1]{Gublerheights}) such that there exists an isomorphism $\varphi\colon L^d\cong f^*L$, where $f\colon X\to X$ is a surjective morphism and $d>1$ is an integer. There exists a collection of metrics $||\cdot||_n$ on $L$ defined as follows, see \cite[\S2.1]{Zhangsmallpoinyts}\[||\cdot||_1=||\cdot||, \quad||\cdot||_n=(\varphi^*f^*||\cdot||_{n-1})^{\frac{1}{d}}\quad(n\ge2).\]
   \begin{theo}\cite[Theorem~2.2]{Zhangsmallpoinyts}\label{thm 1}
       The sequence of metrics $||\cdot||_n$  converges uniformly to a metric $||\cdot||_0$ on $L$, that is, the sequence of functions $\log\frac{||\cdot||_n}{||\cdot||_1}$ converges uniformly  to $\log\frac{||\cdot||_0}{||\cdot||_1}$ on $X(\bar{K})$. Moreover, the metric $||\cdot||_0$ is the unique continuous and bounded metric satisfying the equation $||\cdot||_0=(\varphi^*f^*||\cdot||_0^{\frac{1}{d}})$.  If $\varphi$ is replaced by $\lambda\varphi$ where $\lambda\in K^*$, then the metric $||\cdot||_0$ changes to $|\lambda|^{\frac{1}{d-1}}||\cdot||_0$.  
   \end{theo}
Admissible line bundles on abelian varieties are defined by applying the above theorem to abelian varieties. We firstly construct admissible metrics on symmetric and antisymmetric line bundles as follows.  
   
Let $A$ be an abelian variety over a local field $K$. Fix a symmetric line bundle $L$ on $A$ with a continuous and bounded metric $||\cdot||$.  For any integer $d\ge2$, let $[d]$ be the $d$-multiplication morphism. There is an isomorphism $\varphi\colon L^{\otimes d^2}\cong[d]^*L$ by the theorem of the cube. Theorem~\ref{thm 1} implies the existence of a unique metric $||\cdot||_{d,\varphi}$ on $L$ which satisfies $\varphi^*[d]^*||\cdot||_{d,\varphi}=||\cdot||_{d,\varphi}^{d^2}$. A metric in the set $\pi_d(A,L)=\mathbb{R}^+\cdot||\cdot||_{d,\varphi}$ is called a \textit{$d$-admissible metric}. Note that this set is independent of $\varphi$ by Theorem~\ref{thm 1}. By \cite[Lemma~3.4]{Heinz} and \cite[Remarques~4.8]{Moretmetriquespermises}, the set of $d$-admissible line bundles is independent of $d\ge2$; hence we can denote $\pi(A,L)\coloneqq\pi_d(A,L)$ for any $d$. An element of $\pi(A,L)$ is called an \textit{admissible metric} on $L$.
       

For an antisymmetric line bundle $L^{'}$ on $A$, there is an isomorphism $\varphi\colon L^{'\otimes d}\cong[d]^*L^{'}$. Admissible metrics on $L^{'}$ are defined by $\pi(A,L^{'})=\mathbb{R}^+\cdot||\cdot||_{d,\varphi}$, where $||\cdot||_{d,\varphi}$ is the unique metric that satisfies $||\cdot||_{d,\varphi}^d=\varphi^*[d]^*||\cdot||_{d,\varphi}$. Admissible metrics on any line bundles over an abelian variety are defined as follows.
\begin{defi}
Let $L$ be a line bundle on an abelian variety $A$ over a local field $K$. There is a decomposition $L^{\otimes2}\cong L_1\otimes L_2$, where $L_1=L\otimes[-1]^*L$ is symmetric and $L_2=L\otimes[-1]^*L^{\otimes-1}$ is antisymmetric. A metric $||\cdot||$ on $L$ is called \textit{admissible} if there is an admissible metric $||\cdot||_i$ on $L_i$ for $i\in\{1,2\}$ such that $||\cdot||^2=||\cdot||_1\otimes||\cdot||_2$.     
\end{defi}
Furthermore, we define admissible metrized line bundles on a curve over a local field $K$ by pulling back admissible metrized line bundles from the Jacobian of the curve. Let $C/K$ be a curve of genus $g\ge 1$, and let $J(C)$ be its Jacobian variety. Assume there exists a point $\alpha\in C(K)$, otherwise, we can extend the field $K$. There is an embedding $j_\alpha\colon C\hookrightarrow J(C)$ defined over $K$ which maps a point $x\in C(K)$ to $[x-\alpha]$. 

Let $L$ be a line bundle on $C$ of degree $d$ and let $\theta=\mathcal{O}_J(\Theta^\alpha)$. Since $\deg(j_\alpha^*\theta)=g$, the line bundle $L^{\otimes g}\otimes(j_\alpha^*\theta^{\otimes d})^{-1}$ has degree zero. There exists $N\in\operatorname{Pic}^0(J(C))$ such that $j_\alpha^*N\cong L^{\otimes g}\otimes(j_\alpha^*\theta^{\otimes d})^{-1}$ because $j_\alpha^*\colon\operatorname{Pic}^0(J(C))\to\operatorname{Pic}^0(C)$ is an isomorphism. Setting $M_g=\theta^{\otimes d}\otimes N$, we have $j_\alpha^*M_g\cong L^{\otimes g}$, and the class of $M_g$ in $\operatorname{NS}(J(C))$ is $d[\theta]$.
\begin{defi}
    The set of admissible metrics $\pi(C,L)$ for any line bundle $L$ on $C$ is defined as\[\pi(C,L)\coloneqq\{ ||\cdot||:||\cdot||^{\otimes g}\in j_\alpha^*\pi(J(C),M_g)\}.\] It is shown in \cite[\S4]{Heinz} that $\pi(C,L)$ is independent of the choice of the point $\alpha$.
\end{defi}
We now define the bi-admissible line bundles on $C\times C$,  which are induced from admissible line bundles on $C$.  
\begin{defi}\cite[Definition~4.7]{Moretmetriquespermises},\cite[Definition~4.2]{Heinz}
Let $C$ be a  curve of genus $g\ge 1$ over a local field $K$, and let $L$ be a line bundle on $C\times C$. A metric on $L$ is called \textit{bi-admissible} if for any embedding $j\colon C\hookrightarrow C\times C$ such that one of the compositions $\operatorname{pr}_i\circ j$ is the identity and the other one is constant, the metric $j^*||\cdot||$ is admissible on $i^*L$. 
\end{defi}
A bi-admissible line bundle can be constructed by pulling back an admissible line bundle on $C$ via a projection morphism as follows.  

\begin{theo}\cite[Remark~4.8]{Moretmetriquespermises},\cite[Section~4]{Heinz}
Let $L$ be a metrized line bundle on $C$ with an admissible metric $||\cdot||$. Then $\operatorname{pr}_i^*||\cdot||$ is a bi-admissible metric on $\operatorname{pr}_i^*L$, where $\operatorname{pr}_i\colon C\times C\to C$ is the  $i$th projection morphism.
\end{theo}

In general, there exists a unique bi-admissible metric up to a constant factor on any line bundle $L$ on $C\times C$, see \cite[Lemma~2.2.1, Lemma~2.2.2, Lemma~2.2.3]{ZhangGrossSchoen}. We consider bi-admissible metrics on the diagonal bundle $\mathcal{O}_{C^2}(\Delta)$ as follows.

\begin{theo}\cite[Proposition~4.3]{Heinz}\label{biadmissible}
        The diagonal bundle $\mathcal{O}_{C^2}(\Delta)$ admits a bi-admissible symmetric metric $||\cdot||$. Moreover, $\Delta^*||\cdot||$ is an admissible metric on $\Delta^*\mathcal{O}_{C^2}(\Delta)\cong\omega^{-1}_C$ via the diagonal embedding $\Delta\colon C\to C\times C,x\mapsto(x,x)$.
    \end{theo}
The admissible pairing on $C$ can be constructed using the bi-admissible metrics on $\mathcal{O}_{C^2}(\Delta)$ as follows.
    \begin{theo}\cite[Theorem~4.4]{Heinz},\cite[\S4.7]{Zhangadmissiblepairing}
        Let $||\cdot||$ be a bi-admissible metric on $\mathcal{O}_{C^2}(\Delta)$, then the pairing \[\langle\,,\,\rangle\colon\{(D_1,D_2)\in\operatorname{Div}\bar{C}\times\operatorname{Div}\bar{C}\colon|D_1|\cap|D_2|=\emptyset\}\longrightarrow\mathbb{R}
        \]defined by\[\left(\sum_in_ix_i,\sum_jm_jy_j\right)\mapsto-\sum_{i,j}n_im_j\log||1_\Delta(x_i,y_j)||\] is equal to the admissible pairing. 
    \end{theo}
The following lemma is used in the proof of Theorem~\ref{biadmissible}; it will be applied in \S\ref{ch6}. 
\begin{lemma}\label{lemmaHeinz}\cite[Lemma~4.5]{Heinz}
Let $L$ be a line bundle on $\bar{C}$, and let $s$ be a nonzero rational section of $L$ defined over a Zariski open subset $U\subset\bar{C}$. If $L$ has an admissible metric $||\cdot||$, then there is a constant $c\in\mathbb{R}$ such that \[\langle\operatorname{div}(s),x\rangle_{\operatorname{ad}}=-\log||s(x)||+c\] for any $x\in U\setminus|\operatorname{div}(s)|$.
\end{lemma}
\subsection{Explicit construction of the admissible pairing}\label{3.2}

Let $C$ be a curve of positive genus $g$ over a non-archimedean local field $K$. The admissible pairing on $C$ is constructed in terms of intersection theory on a proper regular model of $C$ and harmonic analysis on the reduction graph of $C$, see \cite[Theorem~4.6]{Zhangadmissiblepairing}. The details of the construction are given below.   



By the semistable reduction theorem, there always exists a finite extension $E$ of $K$ in a fixed algebraic closure $\bar{K}$ such that $C$ has semistable reduction over $E$. Let $\mathcal{C}_E$ be the minimal regular model of $C$ over $\operatorname{Spec }\mathcal{O}_E$; then its special fiber $\mathcal{C}_{E,s}$ is a semistable curve. For any point $P\in C(E)$, let $\bar{P}$ be its Zariski closure in $\mathcal{C}_E$. The absolute intersection multiplicity $(P\cdot Q)$ for distinct points $P,Q\in C(E)$ is defined as 
\begin{equation}
    (P\cdot Q)\coloneq i_E(\bar{P},\bar{Q})/[E:K]
\end{equation} where $i_E(\bar{P},\bar{Q})$ is the intersection multiplicity on $\mathcal{C}_E$. It is easy to verify that $(\,\cdot\,)$ is independent of the unramified field extension. Hence, the intersection multiplicity $(P\cdot Q)$ is well defined for any distinct points $P,Q\in C(\bar{K})$. 

Now, let $G_C$ be the metrized graph given by the dual graph of the special fiber of $\mathcal{C}_E$ with length $\frac{1}{[E:K]}$ on each edge, see \cite[\S2]{Chinburgcapacity}. To be more precise, the set of vertices $V(G_C)$ consists of all irreducible components of the special fiber. Each edge in the set of edges $E(G_C)$ corresponds to an intersection point between two irreducible components. It is well known that $G_C$ does not depend on the field extension up to equivalence of metrized graphs. The specialization of $P\in C(\bar{K})$ is denoted by $\operatorname{sp}(P)\in G_C$.

For any vertex $\Gamma\in V(G_C)$, we define $a_\Gamma\coloneqq-\Gamma^2+2p_a(\Gamma)-2$ where $\Gamma^2$ is the self-intersection number of the associated irreducible component and $p_a(\Gamma)$ is its arithmetic genus. The \textit{canonical divisor} on $G_C$ is defined by \[\label{canonical}K_G\coloneqq\sum a_\Gamma\Gamma,\] where the sum is over all $\Gamma$ in $V(G_C)$. It is easy to prove that $\deg(K_G)=2g-2$. 



Applying \cite[Theorem~3.2]{Zhangadmissiblepairing} to the canonical divisor $K_G$, we obtain the Green function $G(x,y)$ associated to $K_G$ on $G_C\times G_C$, which satisfies \[c+G(K_G,x)+G(x,x)=0\] where $c$ is a constant. 

For any distinct points $P,Q\in C(\bar{K})$, an explicit formula for the admissible pairing is given as follows, see \cite[Theorem~4.6]{Zhangadmissiblepairing} 
\begin{equation}\label{admissiblepairing}
\langle P,Q\rangle_{\operatorname{ad}}=(P\cdot Q)+G(\operatorname{sp}(P),\operatorname{sp}(Q)).
\end{equation}

Let $L$ be the negative of the intersection matrix $M$ on the special fiber $\mathcal{C}_s$, and let $L^+$ be the Moore-Penrose pseudo inverse matrix of $L$, see \cite[\S2]{CinkirExplicitGreenfun}. 
We denote the bilinear pairing associated to $L^{+}$ by $L^{+}(D,E)$ for any $D,E\in\mathbb{Q}^{V(G)}$. An explicit formula for $G(\,,\,)$ on $V(G_C)\times V(G_C)$ is given as follows.     
\begin{lemma}\label{greenfun}
There is a constant $\kappa$ such that the following equality holds for any two vertices  $\Gamma_1,\Gamma_2\in V(G)$ 
     \[G(\Gamma_1,\Gamma_2)=\frac{1}{2g}(-L^+(K_G,\Gamma_1+\Gamma_2)-L^+(\Gamma_1,\Gamma_1)-L^+(\Gamma_2,\Gamma_2)+ 2gL^+(\Gamma_1,\Gamma_2))+\kappa.\]
 \end{lemma}
 \begin{proof}
     We regard a connected metrized graph as a resistive electrical network to perform harmonic analysis on it, see \cite{Chinburgcapacity}. This means that there is an effective resistance function  $r(\Gamma_1,\Gamma_2)$ and a voltage function $j_{\Gamma_i}(\Gamma_1,\Gamma_2)$ between two vertices. We have the following identity, see \cite[Equation~7]{CinkirExplicitGreenfun} \[G(\Gamma_1,\Gamma_2)=\frac{1}{2g}\left(\sum_i a_i j_{\Gamma_i}(\Gamma_1,\Gamma_2)-r(\Gamma_1,\Gamma_2)\right) +\kappa_1\] where $\kappa_1$ is a constant. Moreover, the formulas for $j_{\Gamma_i}(\Gamma_1,\Gamma_2)$ and $r(\Gamma_1,\Gamma_2)$ are given by \cite[Lemma~6.1]{CinkirExplicitGreenfun} as follows
     \begin{equation*}
        \begin{aligned}
             r(\Gamma_1,\Gamma_2)&=L^+(\Gamma_1,\Gamma_1)-2L^+(\Gamma_1,\Gamma_2)+L^+(\Gamma_2,\Gamma_2);\\
             j_{\Gamma_i}(\Gamma_1,\Gamma_2)&=L^+(\Gamma_i,\Gamma_i)-L^+(\Gamma_i,\Gamma_1)-L^+(\Gamma_i,\Gamma_2)+L^+(\Gamma_1,\Gamma_2).
     \end{aligned} 
     \end{equation*}
By bilinearity of $L^+$, we have $L^+(K_G,\Gamma^{'})=L^+(\sum_i a_i\Gamma_i,\Gamma^{'})=\sum_ia_i L^+(\Gamma_i,\Gamma^{'})$. Combining the above formulas, we obtain \[\sum_i a_i j_{\Gamma_i}(\Gamma_1,\Gamma_2)=-L^+(K_G,\Gamma_1+\Gamma_2)+\sum_i a_i\left(L^+(\Gamma_i,\Gamma_i)+L^+(\Gamma_1,\Gamma_2)\right)\] where $\sum_ia_i L^+(\Gamma_i,\Gamma_i)$ is a constant. The proof is completed by combining the above explicit equations.
 \end{proof}

\section{Pulling back the Theta line bundle }\label{ch5}
In this section, we determine the pullback of the theta line bundle via two different morphisms. Let $C$ be a curve of genus $g$ over a field $K$ with a $K$-rational point $\alpha$. First, we define the generalized Abel--Jacobi map $\varphi\colon C^g\to J(C)$. Theorem~\ref{theo1n} computes $\varphi^*\theta^\alpha$ using the properties of the determinant of cohomology, see \cite[\S1]{Moretmetriquespermises}. When $C$ is an odd-degree hyperelliptic curve of genus $g$, Theorem~\ref{coro} computes the divisor $\varphi^*\Theta^\infty$ via Theorem~\ref{theo1n}. Furthermore, we define the difference morphism $\psi\colon C\times C\to J(C)$. In particular, when $C$ is a curve of genus $2$, the generalized Abel--Jacobi map $\varphi$ is essentially equal to the difference morphism $\psi$. Hence $\psi^*\theta^\alpha$ can also be computed using Theorem~\ref{theo1n}.                 
\subsection{Pulling back the Theta line bundle via the generalized Abel--Jacobi map }
Let $C$ be a curve of genus $g$ over a field $K$ with a $K$-rational point $\alpha$. The generalized Abel--Jacobi map $\varphi\colon C^g\to J(C)$ is defined by \[(x_1,...,x_g)\mapsto [
\sum_{i=1}^g x_i-g\alpha].\]
We give a formula for the pullback $\varphi^*\theta^\alpha$  and make it explicit for $g=2$.   
\begin{theo}\label{theo1n} We have an isomorphism \[\varphi^*\theta^\alpha\cong\bigotimes_{i=1}^g\operatorname{pr}_i^*(\omega_C\otimes\mathcal{O}_C(\alpha))\otimes\bigotimes_{1\leq i<j\leq g}\operatorname{pr}_{i,j}^*\mathcal{O}_{C^2}(\Delta)^{-1} \] where $\operatorname{pr}_{i,j}\colon C^g\to C^2$ maps $(x_1,...,x_g)$ to $(x_i,x_j)$ and $
\Delta$ is the diagonal divisor on $C^2$.  
\end{theo}
\begin{proof}
Consider the following diagram
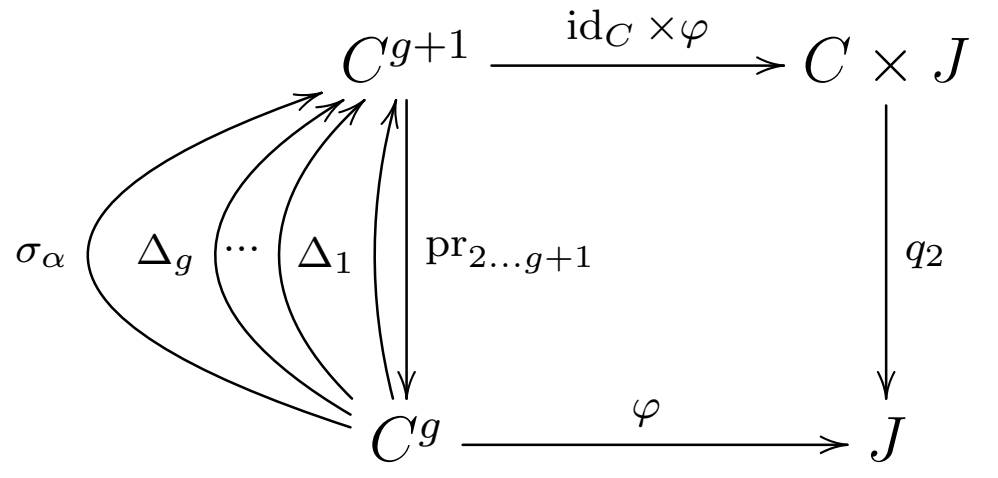
\[\xymatrix{ C^{g+1} \ar[rr]^-{\operatorname{id}_C\times\varphi}\ar[dd]^{\operatorname{pr}_{2...g+1}}&&C\times J\ar[dd]^{q_2}\\\\C^g\ar[rr]^\varphi\ar@/^/[uu]^{\Delta_1}\ar@/^2pc/[uu]^-{...}\ar@/^3pc/[uu]^{\Delta_g}\ar@/^5pc/[uu]^{\sigma_\alpha}&&J}\]
where the morphisms $\Delta_i$ are defined by $C^g\to C^{g+1}\colon(x_1,...,x_g)\to (x_i,x_1,...,x_g)$ for $i\in \{1,\ldots,g\}$ and the morphism $\sigma_\alpha$ is defined by  $C^g\to C^{g+1}\colon(x_1,...,x_g)\to(\alpha,x_1,...,x_g)$ and $q_2$ is the second projection morphism. The notations $\Delta_i$ and $\sigma_\alpha$ also denote the corresponding divisors on $C^{g+1}$ whose support is the image of the morphisms $\Delta_i$ and $\sigma_\alpha$. 

We have the following isomorphism   
\begin{equation}\label{main2}
    \varphi^*\theta^\alpha=\varphi^*(\operatorname{detR}q_{2*}\mathcal{U}^\alpha)^{-1} \cong(\operatorname{detR}\operatorname{pr}_{2...g+1*}(\operatorname{id}\times\varphi)^*\mathcal{U}^\alpha)^{-1},
\end{equation}
where $\mathcal{U}^\alpha$ is the universal line bundle associated to $\alpha$. The universality of $\mathcal{U}^\alpha$ implies the isomorphism   
        \begin{equation}\label{e1}
            (\operatorname{id}\times\varphi)^*\mathcal{U}^\alpha\cong\mathcal{O}_{C^{g+1}}(\Delta_1+...+\Delta_g-g\sigma_\alpha)\otimes \operatorname{pr}_{2...g+1}^*\sigma_\alpha^*\mathcal{O}_{C^{g+1}}(g\sigma_\alpha-\Delta_1-...-\Delta_g).
        \end{equation} 
Moreover, we obtain the following isomorphisms 
\begin{equation}\begin{aligned}\label{e2}
    \sigma_\alpha^*\mathcal{O}_{C^{g+1}}(g\sigma_\alpha-\Delta_1-\cdots-\Delta_g)&\cong c_\alpha^*{\omega_C}^{-g}\otimes \operatorname{pr}_1^{*}\mathcal{O}_C(\alpha)^{-1}\otimes\cdots\otimes  \operatorname{pr}_g^{*}\mathcal{O}_C(\alpha)^{-1}\\
    &\cong\operatorname{pr}_1^{*}\mathcal{O}_C(\alpha)^{-1}\otimes\cdots\otimes \operatorname{pr}_g^{*}\mathcal{O}_C(\alpha)^{-1}\\
\end{aligned}  
\end{equation} 
since $c_\alpha^*{\omega_C}\cong\mathcal{O}_{C^g}$, where $c_\alpha\colon C^g\to C,\,(x_1,...,x_g)\to\alpha$ is a constant map.

We define the line bundles $L_i\coloneqq\mathcal{O}_{C^{g+1}}(\Delta_1+...+\Delta_g-i\sigma_\alpha)$ for any $i\leq g$ on $C^{g+1}$. The following isomorphism is obtained by \eqref{e1}, \eqref{e2}.
\begin{equation}\label{Lg}
    \begin{aligned}
        \operatorname{detR}\operatorname{pr}_{2...g+1*}(\operatorname{id}\times\varphi)^*\mathcal{U}^\alpha &\cong\operatorname{detR}\operatorname{pr}_{2...g+1*}\left(L_g\otimes\operatorname{pr}_{2...g+1}^*\sigma_\alpha^*\mathcal{O}_{C^{g+1}}(g\sigma_\alpha-\Delta_1-...-\Delta_g)\right)\\
        &\cong\operatorname{detR}\operatorname{pr}_{2...g+1*}L_g\otimes\left(\sigma_\alpha^*\mathcal{O}_{C^{g+1}}(g\sigma_\alpha-\Delta_1-...-\Delta_g)\right)^{1-g}\\ 
        &\cong\operatorname{detR}\operatorname{pr}_{2...g+1*}L_g\otimes \operatorname{pr}_1^*\mathcal{O}_C(\alpha)^{g-1}\otimes...\otimes \operatorname{pr}_g^*\mathcal{O}_C(\alpha)^{g-1}.
    \end{aligned}
\end{equation}
Next, we compute the line bundle $\operatorname{detR}\operatorname{pr}_{2...g+1*}L_g$. Consider the exact sequence\begin{equation}\label{eq1}
        0\to\mathcal{I}_{\sigma_\alpha}\to \mathcal{O}_{C^{g+1}}\to \mathcal{O}_{C^{g+1}}/{\mathcal{I}_{\sigma_\alpha}}\to 0
    \end{equation}
where $\mathcal{I}_{\sigma_\alpha}$ is the ideal sheaf of the divisor $\sigma_\alpha$. For the line bundle $L_i$ where $0\leq i\leq g-1$, tensoring $L_i$ with the exact sequence \eqref{eq1} implies the isomorphism
\begin{equation}\label{e3}
    \operatorname{det R}\operatorname{pr}_{2...g+1*}L_{i+1}\cong\operatorname{det R}\operatorname{pr}_{2...g+1*}L_{i}\otimes\operatorname{pr}_1^*\mathcal{O}_C(\alpha)^{-1}\otimes...\otimes \operatorname{pr}_g^*\mathcal{O}_C(\alpha)^{-1}. 
\end{equation} 
The isomorphism $\operatorname{det R}\operatorname{pr}_{2...g+1*}(L_i\otimes \mathcal{O}_{C^{g+1}}/\mathcal{I}_\alpha)\cong\operatorname{pr}_1^*\mathcal{O}_C(\alpha)\otimes...\otimes \operatorname{pr}_g^*\mathcal{O}_C(\alpha)$ for $i>0$ is given by the projection formula. We compute $\operatorname{det R}\operatorname{pr}_{2...g+1*}L_g$ using \eqref{e3} inductively 
\begin{equation}
    \begin{split}
        \operatorname{det R}\operatorname{pr}_{2...g+1*}L_g&\cong\operatorname{det R}\operatorname{pr}_{2...g+1*}L_{g-1}\otimes\operatorname{pr}_1^*\mathcal{O}_C(\alpha)^{-1}\otimes...\otimes \operatorname{pr}_g^*\mathcal{O}_C(\alpha)^{-1}\\
        &\cong\operatorname{det R}\operatorname{pr}_{2...g+1*}L_{g-2}\otimes\operatorname{pr}_1^*\mathcal{O}_C(\alpha)^{-2}\otimes...\otimes \operatorname{pr}_g^*\mathcal{O}_C(\alpha)^{-2}\\
        &\vdotswithin{\cong}\\
        &\cong\operatorname{det R}\operatorname{pr}_{2...g+1*}L_{g-i}\otimes\operatorname{pr}_1^*\mathcal{O}_C(\alpha)^{-i}\otimes...\otimes \operatorname{pr}_g^*\mathcal{O}_C(\alpha)^{-i}\\
         & \vdotswithin{\cong} \\
        &\cong\operatorname{det R}\operatorname{pr}_{2...g+1*}L_0\otimes\operatorname{pr}_1^*\mathcal{O}_C(\alpha)^{-g}\otimes...\otimes\operatorname{pr}_g^*\mathcal{O}_C(\alpha)^{-g} \label{Li}\\
    \end{split}
\end{equation}
Furthermore, the line bundle $\operatorname{det R}p_*L_0$ is computed as follows 
\begin{equation}\label{L0}
    \begin{split}
        \operatorname{det R}\operatorname{pr}_{2...g+1*}\mathcal{O}_{C^{g+1}}(\Delta_1+...+\Delta_g)\cong&\bigotimes_{i=1}^g\Delta_i^*\mathcal{O}_{C^{g+1}}(\sum_{j=1}^i\Delta_j)\\\cong&\bigotimes_{i=1}^g\operatorname{pr}_i^*\omega_C^{-1}\otimes\bigotimes_{1\leq i<j\leq g}\operatorname{pr}_{i,j}^*\mathcal{O}_{C^2}(\Delta).
    \end{split}
\end{equation}
Hence the pullback of the theta line bundle $\theta^\alpha$ is given by 
\begin{equation*}
    \begin{split}
        \varphi^*\theta^\alpha&\cong\varphi^*(\operatorname{det R}q_{2*}\mathcal{U}^\alpha)^{-1}
     \\& \cong(\operatorname{det R}\operatorname{pr}_{2...g+1*}(\operatorname{id}\times\varphi)^*\mathcal{U}^\alpha)^{-1}
    \\&\cong(\operatorname{det R}\operatorname{pr}_{2...g+1*}(L_g\otimes \operatorname{pr}_{2...g+1}^*\sigma_\alpha^*L_g^{-1}))^{-1}\text{ by }\eqref{Lg}
    \\&\cong(\operatorname{det R}\operatorname{pr}_{2...g+1*}\mathcal{O}_{C^{g+1}}(\Delta_1+...+\Delta_g)\otimes\operatorname{pr}_1^*\mathcal{O}_C(\alpha)^{-1}\otimes...\otimes \operatorname{pr}_g^*\mathcal{O}_C(\alpha)^{-1})^{-1}\text{ by }\eqref{Li}
    \\&\cong\bigotimes_{i=1}^g\operatorname{pr}_i^*\left(\omega_C\otimes\mathcal{O}_C(\alpha)\right)\otimes\bigotimes_{1\leq i<j\leq g}\operatorname{pr}_{i,j}^*\mathcal{O}_{C^2}(\Delta)^{-1}\text{ by } \eqref{L0}.          
    \end{split}
\end{equation*}
This completes the proof.
\end{proof}
Let $C$ be a hyperelliptic curve of genus $g$ over a field $K$ and choose the rational point $\alpha$ on $C$ to be the point at infinity $\infty$. The morphism $\varphi$ is given by $(x_1,...,x_g)\mapsto [
\sum_{i=1}^g x_i-g\infty]$. The following formula for the pullback divisor $\varphi^*\Theta^\infty$ is given by applying Theorem~\ref{theo1n} to $C/K$.  
\begin{theo}\label{coro}
Let $\Theta^\infty$ be the theta divisor corresponding to the point $\infty$ on $C/K$, the pullback of $\Theta^\infty$ via the morphism $\varphi$ is given by
\begin{equation}\label{coroeq}
    \varphi^*\Theta^\infty=\sum_{1\leq i<j\leq g}\operatorname{pr}_{i,j}^*(\Delta^{'})+\sum_{i=1}^g\operatorname{pr}_i^*\infty
\end{equation} 
where $\Delta^{'}\coloneqq\{(x,\iota(x))|x\in C\}\subset C\times C$. 
    \end{theo}
\begin{proof}
By the seesaw theorem, there is an isomorphism $\operatorname{pr}_1^*\mathcal{O}_C(2\infty)\otimes\operatorname{pr}_2^*\mathcal{O}_C(2\infty)\cong\mathcal{O}_{C^2}(\Delta+\Delta^{'})$. Let $D$ be an effective divisor of degree $g-1$ on $C$. It is easy to see that $D+\iota D$ is a canonical divisor, where $\iota$ is the hyperelliptic involution. Hence we have an isomorphism  $\omega_C\cong\mathcal{O}_C((2g-2)\infty)$ by letting $D=(g-1)\infty$. Theorem \ref{theo1n} implies the following isomorphism  
\[\varphi^*\theta^\infty\cong\bigotimes_{1\leq i<j\leq g}\operatorname{pr}_{i,j}^*\mathcal{O}_{C^2}(\Delta^{'})\otimes\bigotimes_{i=1}^g\operatorname{pr}_i^*\mathcal{O}_C(\infty).\] 

It remains to compare the support of the divisors in \eqref{coroeq}. For any point $(x_1,..x_i,..\iota x_i,..x_g)\in|\operatorname{pr}_{i,j}^*(\Delta^{'})|$, we have $\varphi(x_1,..x_i,..\iota x_i,..x_g)=[\sum_{j=1,j\ne i}^gx_j-(g-2)\infty]\in|\Theta^\infty|$. For any point $(x_1,...\infty,...x_g)\in|\operatorname{pr}_i^*(\infty)|$, we have $\varphi(x_1,..\infty,..x_g)=[\sum_{j=1,j\ne i}^gx_j-(g-1)\infty]\in|\Theta^\infty|$. Hence $|\sum_{1\leq i<j\leq g}\operatorname{pr}_{i,j}^*(\Delta^{'})+\sum_{1\leq i\leq g}\operatorname{pr}_i^*\infty|\subset|\varphi^*\Theta^\infty|$. 

Then $\varphi^*\Theta^\infty-\sum_{1\leq i<j\leq g}\operatorname{pr}_{i,j}^*(\Delta^{'})+\sum_{1\leq i\leq g}\operatorname{pr}_i^*\infty$ $\varphi^*\Theta^\infty$ is effective and linearly equivalent to $0$. We conclude that $\varphi^*\Theta^\infty=\sum_{1\leq i<j\leq g}\operatorname{pr}_{i,j}^*(\Delta^{'})+\sum_{i=1}^g\operatorname{pr}_i^*\infty$.  
\end{proof}
\subsection{The pullback of the Theta line bundle via the difference morphism}
For a curve $C$ over a field $K$, the \textit{difference morphism} $\psi\colon C\times C\to J(C)$ is defined by\[(x,y)\mapsto[x-y].\]When $C$ is a complex curve, Moret-Bailly obtained a formula for the pullback  $\psi^{*}\theta^\alpha$ with respect to $\alpha\in C(K)$ in \cite[\S4.10.2]{Moretmetriquespermises}, but the details are omitted. Moreover, \cite[Lemma~2.1]{Heinz} generalizes the formula when $C$ is a curve over a noetherian scheme but only states that the proof is analogous to that in \cite{Moretmetriquespermises}. In this section, we give a detailed proof when $C$ is a curve over a field $K$. The proof is similar to the proof of Theorem~\ref{theo1n}, which was in fact inspired by Moret-Bailly's proof.  
\begin{theo}\label{g2thm1}
Let $C$ be a curve of genus $g$ over a field $K$, and let $\alpha\in C(K)$ be a rational point. We have the following isomorphism \[\psi^*\theta^\alpha\cong \mathcal{O}_{C\times C}(\Delta)\otimes\operatorname{pr}_1^*\mathcal{O}_C(\alpha)^{1-g}\otimes\operatorname{pr}_2^*\mathcal{O}_C(\alpha)^{g-1}\otimes\operatorname{pr}_1^*\omega_C.\]
 \end{theo}
\begin{proof}
Consider the diagram
      
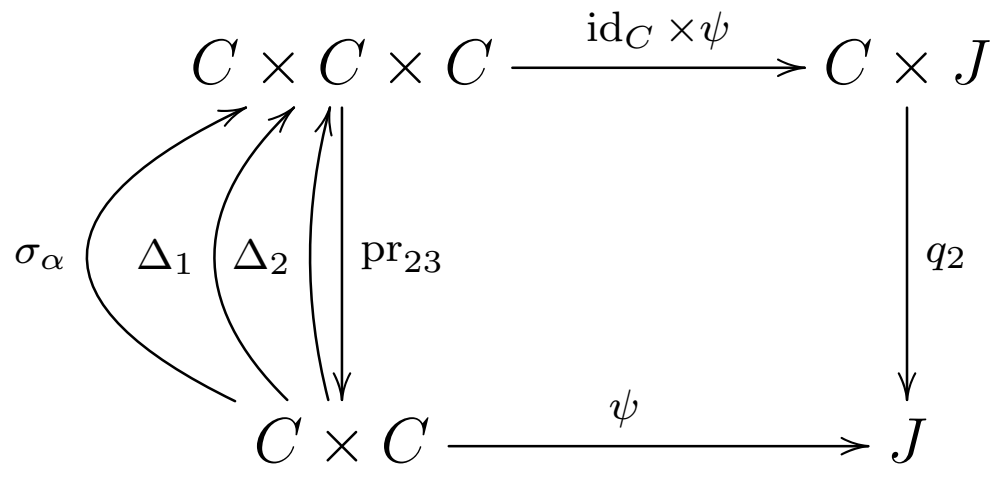
\[\xymatrix{ C\times C\times C  \ar[rr]^-{\operatorname{id}_C\times\psi}\ar[dd]^{\operatorname{pr}_{23}}&&C\times J\ar[dd]^{q_2}\\\\C\times C\ar[rr]^\psi\ar@/^/[uu]^{\Delta_2}\ar@/^2pc/[uu]^{\Delta_1}\ar@/^4pc/[uu]^{\sigma_\alpha}&&J}\]  where $\Delta_1(x,y)\coloneqq(x,x,y)$ , $\Delta_2(x,y)\coloneqq(y,x,y)$, $\sigma_\alpha(x,y)\coloneqq(\alpha,x,y)$. The pullback of the theta line bundle via the morphism $\psi$ is
\begin{equation}\label{iso1}
\begin{aligned}
    \psi^*\theta^\alpha&\cong\psi^*(\operatorname{detR}q_{2*}\mathcal{U}^\alpha)^{-1}\text{ by the definition of }\theta^\alpha\\
    &\cong(\operatorname{detR}\operatorname{pr}_{23*}(\operatorname{id}_C\times\psi)^*\mathcal{U}^\alpha)^{-1}.
\end{aligned}
\end{equation}
Furthermore, the following isomorphism for $(\operatorname{id}_C\times\psi)^*\mathcal{U}^\alpha$ is due to the universality of $\mathcal{U}^\alpha$
\begin{equation*}
    \begin{aligned}
    (\operatorname{id}_C\times\psi)^*\mathcal{U}^\alpha&\cong\mathcal{O}_{C^3}(\Delta_1-\Delta_2)\otimes\operatorname{pr}_{23}^*\sigma_\alpha^*\mathcal{O}_{C^3}(\Delta_2-\Delta_1)\\&\cong\mathcal{O}_{C^3}(\Delta_1-\Delta_2)\otimes\operatorname{pr}_{23}^*(\operatorname{pr}_2^*\mathcal{O}_C(\alpha)\otimes\operatorname{pr}_1^*\mathcal{O}_C(\alpha)^{-1}).    
    \end{aligned}
\end{equation*}
Then we have the following isomorphism
\begin{equation}\label{iso2}
    \begin{aligned}
        \operatorname{detR}\operatorname{pr}_{23*}(\operatorname{id}_C\times\psi)^*\mathcal{U}^\alpha&\cong\operatorname{detR}\operatorname{pr}_{23*}\left(\mathcal{O}_{C^3}(\Delta_1-\Delta_2)\otimes\operatorname{pr}_{23}^*(\operatorname{pr}_2^*\mathcal{O}_C(\alpha)\otimes\operatorname{pr}_1^*\mathcal{O}_C(\alpha)^{-1})\right)\\&\cong\operatorname{detR}\operatorname{pr}_{23*}\mathcal{O}_{C^3}(\Delta_1-\Delta_2)\otimes\left(\operatorname{pr}_2^*\mathcal{O}_C(\alpha)\otimes\operatorname{pr}_1^*\mathcal{O}_C(\alpha)^{-1})\right)^{1-g}.
    \end{aligned}
\end{equation}
We now obtain an isomorphism for $\operatorname{detR}\operatorname{pr}_{23*}\mathcal{O}_{C^3}(\Delta_1-\Delta_2)$ by using the following canonical exact sequence 
\begin{equation}\label{exact}
    0\to\mathcal{I}_{\Delta_2}\to \mathcal{O}_{C^3}\to \mathcal{O}_{\Delta_2}\to 0
\end{equation} 
where $\mathcal{I}_{\Delta_2}= \mathcal{O}_{C^3}(-\Delta_2)$ and $\mathcal{O}_{\Delta_2}= \mathcal{O}_{C^3}/\mathcal{I}_{\Delta_2}$. Tensoring the line bundle $\mathcal{O}_{C^3}(\Delta_1)$ with \eqref{exact}, we have
\begin{equation}\label{g2eq1}
    \operatorname{detR}\operatorname{pr}_{23*}\mathcal{O}_{C^3}(\Delta_1-\Delta_2)\cong \operatorname{detR}\operatorname{pr}_{23*}\mathcal{O}_{C^3}(\Delta_1)\otimes\left(\operatorname{detR}\operatorname{pr}_{23*}(\mathcal{O}_{C^3}/\mathcal{I}_{\Delta_2}\otimes \mathcal{O}_{C^3}(\Delta_1))\right)^{-1}.
\end{equation}    
Here we have the isomorphism $\operatorname{detR}\operatorname{pr}_{23*}\mathcal{O}_{C^3}(\Delta_1)\cong\operatorname{detR}\operatorname{pr}_{23*}\mathcal{O}_{ C^3}\otimes \Delta_1^*\mathcal{O}_{C^3}(\Delta_1)$. The isomorphism  $\Delta_{2*}\mathcal{O}_{\Delta_2}\otimes_{\mathcal{O}_{C^3}} \mathcal{O}_{C^3}(\Delta_1)\cong \Delta_{2*}(\mathcal{O}_{\Delta_2}\otimes_{\mathcal{O}_{\Delta_2}}\Delta^*_2
 \mathcal{O}_{C^3}(\Delta_1))$ is obtained by the projection formula. Combining these isomorphisms yields   
   \begin{equation}\label{iso3} 
     \operatorname{detR}\operatorname{pr}_{23*}(\mathcal{O}_{C^3}/\mathcal{I}_{\Delta_2}\otimes_{\mathcal{O}_{C^3}} \mathcal{O}_{C^3}(\Delta_1))\cong\Delta_2^*\mathcal{O}_{C^3}(\Delta_1). 
    \end{equation}   
Therefore, the following isomorphism for $\operatorname{detR}\operatorname{pr}_{23*}\mathcal{O}_{C^3}(\Delta_1-\Delta_2)$ follows from \eqref{g2eq1} and \eqref{iso3} 
\begin{equation}\label{iso4}
    \begin{aligned}
        \operatorname{detR}\operatorname{pr}_{23*}\mathcal{O}_{C^3}(\Delta_1-\Delta_2)&\cong\operatorname{detR}\operatorname{pr}_{23*}\mathcal{O}_{C^3}\otimes \Delta_1^*\mathcal{O}_{C^3}(\Delta_1)\otimes \Delta_2^*\mathcal{O}_{C^3}(\Delta_1)^{-1}\\&\cong\Delta_1^*\mathcal{O}_{C^3}(\Delta_1)\otimes \Delta_2^*\mathcal{O}_{C^3}(\Delta_1)^{-1}.
    \end{aligned}
\end{equation}  
The pullback of the theta line bundle is given by \eqref{iso1},\eqref{iso2} and \eqref{iso4}.  
\end{proof}

Let $C/K$ be a curve of genus $2$. The difference morphism $\psi$ is equal to the morphism $\varphi$ up to the hyperelliptic involution. We reprove Theorem~\ref{g2thm1} in the case of genus $2$ as follows by using Theorem~\ref{theo1n}.  
\begin{ex}\label{ex 5.4}
    Let $C/K$ be a curve of genus $2$. For any points $p,q\in C$, the divisor $(p-q)$ is linearly equivalent to the divisor $(p+\iota q-2\infty)$ since $\operatorname{div}(x-x(q))=(q+\iota q-2\infty)$, where $\iota$ is the hyperelliptic involution. Hence the difference morphism $\psi$ agrees with the morphism $ \varphi^{'}\colon C\times C\to J,\,(p,q)\mapsto[p+\iota q-2\infty]$. We obtain the following isomorphisms for the pullback $\psi^{*}\theta^\infty$ using the same method as in the proof of Theorem~\ref{theo1n}.     
\begin{equation*}
        \begin{aligned}
             \varphi^{'*}\theta^\infty&\cong\varphi^{'*}(\operatorname{detR}q_{2*}\mathcal{U}^\infty)^{-1}
            \\&\cong(\operatorname{detR}\operatorname{pr}_{23*}(\operatorname{id}\times  \varphi^{'})^*\mathcal{U}^\infty)^{-1}
            \\&\cong(\operatorname{detR}\operatorname{pr}_{23*}(L_2^{'}\otimes\operatorname{pr}_{23}^*\sigma_\infty^*(L_2^{'})^{-1}))^{-1}
            \\&\cong(\operatorname{detR}\operatorname{pr}_{23*}L_2^{'}\otimes c_\infty^*{\omega_C}^{-1}\otimes \operatorname{pr}_1^*\mathcal{O}_C(-\infty)\otimes \operatorname{pr}_2^*\mathcal{O}_C(-\infty))^{-1}
            \\&\cong(\operatorname{detR}\operatorname{pr}_{23*}\mathcal{O}_{C^3}(\Delta_1+\Delta_2^{'})\otimes\operatorname{pr}_1^*\mathcal{O}_C(-\infty)\otimes \operatorname{pr}_2^*\mathcal{O}_C(-\infty))^{-1}
            \\&\cong\operatorname{pr}_1^*(\omega_C\otimes\mathcal{O}_C(\infty))\otimes\operatorname{pr}_2^*(\omega_C\otimes\mathcal{O}_C(\infty))\otimes\mathcal{O}_{C^2}(\Delta^{'})^{-1}
        \end{aligned}
\end{equation*}
    where $\Delta_2^{'}\coloneqq\{(\iota y,x,y)\,|\,(x,y)\in C^2\}$, $\,\Delta^{'}\coloneqq\{(x,\iota x)|x\in C\}$, and $\,L_2^{'}\coloneqq\mathcal{O}_{C^3}(\Delta_1+\Delta_2^{'}-2\sigma_\infty)$. 
    Moreover, there is an isomorphism $\operatorname{pr}_1^*\mathcal{O}_C(2\infty)\otimes\operatorname{pr}_2^*\mathcal{O}_C(2\infty)\cong\mathcal{O}_{C^2}(\Delta+\Delta^{'})$ by the seesaw principle. The following isomorphism is the same as \cite[Lemma~4.10.2]{Moretmetriquespermises} since $\mathcal{O}_C(2\infty)\cong\omega_C$ for a curve of genus $2$ \[ \varphi^{'*}\theta^\infty\cong\mathcal{O}_{C^2}(\Delta)\otimes\operatorname{pr}_1^*\mathcal{O}_C(\infty)\otimes\operatorname{pr}_2^*\mathcal{O}_C(\infty).\]
    \end{ex}
\section{N\'eron functions on Jacobians of odd degree hyperelliptic curves} \label{ch6}
There are several approaches to construct explicit formulas for N\'eron functions on Jacobians of curves of low genus, for instance, see \cite{Flynncanonical}, \cite{Silvermancomputing}, \cite{Stollgenus3}. Let $C/K$ be a hyperelliptic curve of genus $g$ over a local field $K$ of characteristic different from $2$ in this section. We define the Abel--Jacobi embedding with respect to the base point $\infty\in C(K)$ and fix the notations $\Theta\coloneq\Theta^\infty$, $\theta\coloneq\theta^\infty$. The formula for the N\'eron function $\lambda_{\Theta}$ is given by the admissible pairings on $C$ in Theorem~\ref{Thm6.2}. The proof of Theorem~\ref{Thm6.2} is based on the formula for the pullback of the theta divisor in Theorem~\ref{coro} and the relation of the admissible metrized line bundles to the associated N\'eron functions on $J$, see Theorem~\ref{Heinz}. The formula in Theorem~\ref{Thm6.2} makes it possible to compute N\'eron functions explicitly. Furthermore, Theorem~\ref{main} relates Theorem~\ref{Thm6.2} to local N\'eron pairings in the non-archimedean case.        

The following result expresses N\'eron functions in terms of admissible metrics and relates N\'eron functions to admissible pairings.  
\begin{theo}\label{Heinz}
Let $A$ be an abelian variety over a local field $K$, and let $L$ be a line bundle on $A$. Let $||\cdot||$ be an admissible metric on $L$, and let $l$ be a nonzero rational section of $L$. The map $x\mapsto-\log||l(x)||$ is a N\'eron function associated to the divisor $\operatorname{div}(l)$.
\end{theo}
\begin{proof}
This was proved in \cite[Theorem~3.7]{Heinz} when $K$ is a non-archimedean local field. By applying \cite[Remark~3.10]{Moretmetriquespermises}, the proof of \cite[Theorem~3.7]{Heinz} also applies in the archimedean case.    
\end{proof}
It follows from \cite[Theorem~6.1]{MenezesElementary} and the Riemann–Roch theorem that every point $a\in J(C)(K)\setminus|\Theta|$ has a representation $a=[D_a]$ where $D_a\coloneqq Q_1+\cdots+Q_g-g\infty\in\operatorname{Div}^0(C)$ with $Q_i\in C(\overline{K})\setminus\{\infty\}$ and $\iota(Q_i)\ne Q_j$ for all $i\ne j$. Then we have the following formula for the N\'eron function. 
\begin{theo}\label{Thm6.2}
For $a=[D_a]\in J(C)(K)\setminus|\Theta|$, a N\'eron function associated to the divisor $\Theta$ is given by \[\lambda_{\Theta}(a)=\sum_{1\leq i<j\leq g}\langle Q_i,\iota Q_j\rangle_{\operatorname{ad}}+\sum_{i=1}^g\langle Q_i,\infty\rangle_{\operatorname{ad}}.\] 
\end{theo}
\begin{proof}
  We define a morphism $\operatorname{id}\times\iota$ as follows \[\operatorname{id}\times\iota\colon C\times C\to C\times C\quad(x,y)\mapsto(x,\iota y).\]
  It is easy to see that $(\operatorname{id}\times\iota)^*\Delta^{'}=\Delta$ and $(\operatorname{id}\times\iota)^*\Delta=\Delta^{'}$. Theorem~\ref{biadmissible} shows that $\mathcal{O}_{C^2}(\Delta)$ has a bi-admissible metric $||\cdot||$. Requiring $\operatorname{id}\times\iota$ to be an isometry, $\mathcal{O}_{C^2}(\Delta^{'})$ is a metrized line bundle with a bi-admissible metric $||\cdot||^{'}$.

Consider the morphism $\varphi\colon C^g\to J$ defined by $(Q_1,...,Q_g)\mapsto[\sum Q_i-g\infty]$. We have the following equations for a N\'eron function $\lambda_{\Theta}(a)$ according to Theorem~\ref{Heinz}
\begin{equation*}
    \begin{aligned}
        \lambda_{\Theta}(a)&=\lambda_{\Theta}([Q_1+...+Q_g-g\infty])\\
        &=-\log||1_{\Theta}(Q_1+...+Q_g-g\infty)||_{\theta}\\
         &=-\log ||1_{\Theta}(\varphi(Q_1,...,Q_g)||_{\theta}\\
          &=-\log ||\varphi^*1_{\Theta}(Q_1,...,Q_g)||_{\varphi^*\theta}\\
          &=-\sum_{1\leq i<j\leq g}\log ||1_{\operatorname{pr}_{i,j}^*\Delta^{'}}(Q_1,...,Q_g)||_{\operatorname{pr}_{i,j}^*\mathcal{O}_C^2(\Delta^{'})}-\sum_{i=1}^g\log||1_{\operatorname{pr}_i^*\infty}(Q_1,...,Q_g)||_{\operatorname{pr}_i^*\mathcal{O}_C(\infty)}\\
    \end{aligned}
\end{equation*}
where $||\cdot||_{\operatorname{pr}_{i,j}^*\mathcal{O}_C^2(\Delta^{'})}$ is the pullback metric from $||\cdot||^{'}$. Then we have 
\begin{equation*}
    \begin{aligned}
        \log||1_{\operatorname{pr}_{i,j}^*\Delta^{'}}(Q_1,...,Q_g)||_{\operatorname{pr}_{i,j}^*\mathcal{O}_C^2(\Delta^{'})}&=\log ||1_{\Delta^{'}}(Q_i,Q_j)||^{'}\\
        &=\log ||1_{(\operatorname{id}\times\iota)^*\Delta}(Q_i,Q_j)||^{'}\\
        &=\log ||1_\Delta(\operatorname{id}\times\iota(Q_i,Q_j))||_{\mathcal{O}_{C^2}(\Delta)}\\
        &=\log ||1_\Delta(Q_i,\iota Q_j)||_{\mathcal{O}_{C^2}(\Delta)}\\
        &=-\langle Q_i,\iota Q_j\rangle_{\operatorname{ad}}.
    \end{aligned}
\end{equation*}
Finally, we obtain the formula for $\lambda_{\Theta}(a)$ as follows  
\begin{equation*}
    \begin{aligned}
        \lambda_{\Theta}(a)&=\lambda_{\Theta}([Q_1+...+Q_g-g\infty])\\
        &=-\sum_{1\leq i<j\leq g}\log||1_{\operatorname{pr}_{i,j}^*\Delta^{'}}(Q_1,...,Q_g)||_{\operatorname{pr}_{i,j}^*\mathcal{O}_{C^2}(\Delta^{'})}-\sum_{i=1}^g\log||1_{\operatorname{pr}_i^*\infty}(Q_1,...,Q_g)||_{\operatorname{pr}_i^*\mathcal{O}_C(\infty)}\\
        &=\sum_{1\leq i<j\leq g}\langle Q_i,\iota Q_j\rangle_{\operatorname{ad}}+\sum_{i=1}^g\langle Q_i,\infty\rangle_{\operatorname{ad}}.
    \end{aligned}
\end{equation*}
This proves the theorem.
\end{proof}
For a curve $C/K$ of genus $g$, the difference morphism is defined as \[\psi\colon C\times C\to J,\; (x,y)\mapsto[x-y].\] This morphism is surjective if and only if $g\le2$. For curves of genus $2$, an explicit formula for $\lambda_{\Theta}(a)$ is given by \cite[Proposition~6.5]{Bianchicoleman} using the difference morphism. We give another proof of this result using the same method as in the proof of Theorem~\ref{Thm6.2}.  
\begin{coro}\cite[Proposition~6.5]{Bianchicoleman}\label{Bianchi}
Let $C/K$ be an odd-degree curve of genus $2$ over a local field $K$. For $a\in J(C)(K)\setminus|\Theta|$, there exist two distinct points $Q_1,Q_2\in C\setminus\{\infty\}$ such that $a=[Q_1-Q_2]$. Then we have \[\lambda_{\Theta^\infty}(a)=\langle Q_1,Q_2\rangle_{\operatorname{ad}}+\langle Q_1,\infty\rangle_{\operatorname{ad}}+\langle Q_2,\infty\rangle_{\operatorname{ad}}.\]  
\end{coro}
\begin{proof}
For a point $Q\in C$, the equality $-\log||1_\infty(Q)||_{\mathcal{O}_C(\infty)}=\langle Q,\infty\rangle_{\operatorname{ad}}$ follows from Lemma~\ref{lemmaHeinz}. Applying Example~\ref{ex 5.4}, we obtain
\begin{equation*}
    \begin{aligned}
        \lambda_{\Theta}(a)&=\lambda_{\Theta}([Q_1-Q_2])\\
        &=-\log||1_{\Theta}([Q_1-Q_2])||_{\theta}\\
        &=-\log||1_{\Theta}(\psi(Q_1,Q_2)||_{\theta}\\
        &=-\log||\psi^*1_{\Theta}(Q_1,Q_2)||_{\psi^*\theta}\\
        &=-\log||1_\Delta(Q_1,Q_2)||_{\mathcal{O}_{C^2}(\Delta)}-\log||1_{\operatorname{pr}_1^*\infty}(Q_1,Q_2)||_{\operatorname{pr}_1^*\mathcal{O}_C(\infty)}-\log||1_{\operatorname{pr}_2^*\infty}(Q_1,Q_2)||_{\operatorname{pr}_2^*\mathcal{O}_C(\infty)}\\
        &=\langle Q_1,Q_2\rangle_{\operatorname{ad}}-\log||1_\infty(\operatorname{pr}_1(Q_1,Q_2))||_{\mathcal{O}_C(\infty)}-\log||1_\infty(\operatorname{pr}_2(Q_1,Q_2))||_{\mathcal{O}_C(\infty)}\\
        &=\langle Q_1,Q_2\rangle_{\operatorname{ad}}-\log||1_\infty(Q_1)||_{\mathcal{O}_C(\infty)}-\log||1_\infty(Q_2)||_{\mathcal{O}_C(\infty)}\\
        &=\langle Q_1,Q_2\rangle_{\operatorname{ad}}+\langle Q_1,\infty\rangle_{\operatorname{ad}}+\langle Q_2,\infty\rangle_{\operatorname{ad}}.
    \end{aligned}
\end{equation*}
This completes the proof.
\end{proof}
\begin{rem}
The two points $Q_1,Q_2$ are unique up to the hyperelliptic involution on $C$, namely the different expressions for $a$ are $[Q_1-Q_2]$ or $[\iota Q_2-\iota Q_1]$. Substituting the points $\iota Q_1$ and $\iota Q_2$ into the above computation yields the same result because the admissible pairing is $\operatorname{Aut}(C)$-invariant.
\end{rem}
\section{Comparison with local N\'eron pairings on the curve}\label{ch7}
Throughout this section, we retain the notation introduced in \S\ref{ch6} except that $K$ now is a non-archimedean local field of characteristic different from $2$ such that $C$ has semistable reduction over $K$.  In this section, we define a local height function $\lambda_N$ using the extended N\'eron pairings $\langle\,,\,\rangle_{N,t}$ on $C$ with respect to a certain set-theoretical section $t$ of the tangent bundle of $C$. Then Theorem~\ref{main} shows that the N\'eron function $\lambda_{2\Theta}$ with respect to the divisor $2\Theta$ on $J(C)$ is equal to the local height function $\lambda_N$ up to an explicit correction term. Moreover, we introduce Coleman--Gross heights and $p$-adic N\'eron functions. We show that Theorem~\ref{padiccomparision} is a $p$-adic analogue of Theorem~\ref{main} and generalizes \cite[Theorem~6.3]{Bianchicoleman}.            
\begin{defi}
 Let $t$ be a set-theoretic section of the tangent bundle of $C$ determined by duality by the choice of the differential $\frac{dx}{2y}$ at all affine points and $\frac{x^{g-1}dx}{2y}$ at the point at infinity $\infty$. For a point $a=[D_a]\in J(C)(K)\setminus|\Theta|$, we define\[\lambda_N(a)\coloneqq\langle D_a,D_a\rangle_{N,t},\] where $D_a$ is defined in \S\ref{ch6}.    
\end{defi}
The choice of tangent vectors in the above definition follows the conventions in \cite[\S5]{BalakrishnanQuadratic}.
\begin{theo}\label{main}
Let $C$ be an odd-degree hyperelliptic curve with semistable reduction over a non-archimedean local field $K$ of characteristic different from $2$, and let $a=[D_a]\in J(C)(K)\setminus|\Theta|$. We have the following identity
\begin{equation}\label{maineq}
        \lambda_N(a)+\lambda_{2\Theta}(a)\equiv-2\sum_{1\leq i<j\leq g}\log|f_j(Q_i-\infty,t)|
\end{equation}
where $f_j$ is the function $x-x(Q_j)$ and $\equiv$ means it is an identity up to an additive constant which is independent of $a$. The notation $f_j(Q_i-\infty,t)$ denotes the extended value defined in \S\ref{ch2.2}.      
\end{theo}

\begin{rem}
Instead of the representation $a=[Q_1+Q_2-2\infty]$, the representation $a=[Q_1-Q_2]$ is used in \cite[\S6]{Bianchicoleman} to express points on $J(C)(K)\setminus|\Theta|$ when $C$ is a curve of genus $2$. 
\end{rem}
\subsection{The proof of Theorem~\ref{main}}
Let $G_C$ be the metrized graph given by the dual graph of the special fiber of $\mathcal{C}_K$ and let $K_G$ be the canonical divisor on $G_C$. Let $L$ be the negative of the intersection matrix $M$ on the special fiber $\mathcal{C}_{K,s}$, and let $L^+$ be the Moore-Penrose pseudo inverse matrix of $L$. The following lemmas give formulas for arithmetic intersection numbers using the $L^+$ pairing on hyperelliptic curves, which will be used frequently in the proof of Theorem~\ref{main}. 
\begin{lemma}\label{lemma difference}
Let $C/K$ be a hyperelliptic curve of genus $g$. For any point $Q\in C(K)\setminus\{\infty\}$, we have
\[((\bar{Q}-\iota\bar{Q})\cdot\bar{\infty})
             =\frac{1}{2g}(L^+(K_G,\Gamma_{Q}-\Gamma_{\iota Q})-L^+(\Gamma_{\iota Q},\Gamma_{\iota Q})+L^+(\Gamma_Q,\Gamma_Q)+ 2gL^+(\Gamma_{\iota Q}-\Gamma_Q,\Gamma_\infty)).\]    
 \end{lemma}
\begin{proof}
The equality $\langle Q,\infty\rangle_{\operatorname{ad}}=\langle \iota Q,\infty\rangle_{\operatorname{ad}}$ holds because the admissible pairing is $\operatorname{Aut}(C)$-invariant on $\operatorname{Div}(C)$. By Zhang's construction of the admissible pairing, see \S\ref{3.2}, we have $\langle Q,\infty\rangle_{\operatorname{ad}}=(\bar{Q}\cdot\bar{\infty})+G(\Gamma_Q,\Gamma_\infty)$
and $\langle \iota Q,\infty\rangle_{\operatorname{ad}}=(\iota\bar{Q}\cdot\bar{\infty})+G(\Gamma_{\iota Q},\Gamma_\infty)$. Comparing these two admissible pairings yields $(\bar{Q}-\iota\bar{ Q}\cdot\bar{\infty})=G(\Gamma_{\iota Q},\Gamma_\infty)-G(\Gamma_Q,\Gamma_\infty)$. Lemma~\ref{greenfun} implies the following formula       
\begin{equation*}
    \begin{aligned}
((\bar{Q}-\iota\bar{Q})\cdot\bar{\infty})&=G(\Gamma_{\iota Q},\Gamma_\infty)-G(\Gamma_Q,\Gamma_\infty)\\
        &=\frac{1}{2g}(L^+(K_G,\Gamma_Q-\Gamma_{\iota Q})-L^+(\Gamma_{\iota Q},\Gamma_{\iota Q})+L^+(\Gamma_Q,\Gamma_Q)+ 2gL^+(\Gamma_{\iota Q}-\Gamma_Q,\Gamma_\infty)).\\
    \end{aligned}
\end{equation*} 
This completes the proof.
\end{proof}
\begin{rem}
 The involution morphism $\iota$ on $C$ extends to an automorphism on the semistable model of $C$. Hence, we have $((\bar{Q}-\iota\bar{Q})\cdot\bar{\infty})=0$ but we will not use this in the following. 
\end{rem}
Let $D$ be a divisor of degree zero on $C$. Proposition~\ref{othogonal} states that there is a $\mathbb{Q}$-divisor $\Phi(D)$ on the proper regular semistable model $\mathcal{C}$ such that $\bar{D}+\Phi(D)$ has trivial intersections with all vertical divisors, where $\bar{D}$ is the Zariski closure of $D$ in $\mathcal{C}$.
\begin{lemma}\label{lemma L}\cite[Proposition~7.4]{HolmesAsymptotics}
    Let $D=\sum_in_iP_i$ and $E=\sum_im_iQ_i$ be divisors in $\operatorname{Div}^0(C(K))$. Let $\Gamma_{P_i}$ and $\Gamma_{Q_i}$ be the specializations of $P_i$ and $Q_i$ on the dual graph $G(C)$. Denote $\sum_in_i\Gamma_{P_i}$ by $\Gamma_D$ and $\sum_im_i\Gamma_{Q_i}$ by $\Gamma_E$. We have an equality\[(\Phi(D)\cdot\Phi(E))=-L^+(\Gamma_D,\Gamma_E).\]
\end{lemma}
To prove Theorem~\ref{main}, we need to give explicit formulas for $\lambda_N(a)$ and $\lambda_{2\Theta}(a)$ in terms of intersection multiplicities and the $L^+$ pairings. They are given in the following lemmas.   
\begin{lemma}
An explicit formula for $\lambda_N(a)$ is given by
    \begin{equation}\label{lambda}
        \begin{aligned}
          \lambda_N(a)&\equiv(2-4g)(\sum_{i=1}^g\bar{Q}_i\cdot\bar{\infty})+2\sum_{1\leq i<j\leq g}(\bar{Q_i}\cdot\bar{Q_j})_t+\sum_{i=1}^gL^+(K_G,\Gamma_{Q_i})+(2-4g)L^+(\sum_{i=1}^g\Gamma_{Q_i},\Gamma_\infty)\\ 
          &\quad+\sum_{i=1}^gL^+(\Gamma_{Q_i},\Gamma_{Q_i})+2\sum_{1\leq i<j\leq g}L^+(\Gamma_{Q_i},\Gamma_{Q_j}).
        \end{aligned}
    \end{equation}    
\end{lemma}
\begin{proof}
By the definition of $\lambda_N$ and Lemma~\ref{lemma L}, we have 
      \begin{equation*}
          \begin{aligned}
            &\quad\lambda_N(a)\\
            &=(\bar{Q}_1+...+\bar{Q}_g-g\bar{\infty}\cdot\bar{Q}_1+...+\bar{Q}_g-g\bar{\infty})_t-(\Phi(Q_1+...+Q_g-g\infty)\cdot\Phi(Q_1+...+Q_g-g\infty))\\            &=\sum_{i=1}^g(\bar{Q}_i\cdot\bar{Q}_i)_{t}+g^2(\bar{\infty}\cdot\bar{\infty})_{t}+2\sum_{1\leq i<j\leq g}(\bar{Q_i}\cdot\bar{Q_j})_t-2g\sum_{i=1}^g(\bar{Q}_i\cdot\bar{\infty})\\
            &\quad+L^+(\sum_{i=1}^g\Gamma_{Q_i}-g\Gamma_\infty,\sum_{i=1}^g\Gamma_{Q_i}-g\Gamma_\infty).\\
          \end{aligned}
      \end{equation*}
We express $(\bar{Q}_i\cdot\bar{Q}_i)_{t}$ in terms of $L^+$ and intersection multiplicity as follows. The divisor of the differential form $\omega=\frac{dx}{2y}$ on the regular semistable model $\mathcal{C}$ is $\operatorname{div}(\omega)=(2g-2)\bar{\infty}+V$ where $V$ is a vertical divisor on $\mathcal{C}$. It follows from Lemma~\ref{adjunction} that 
 \begin{equation*}
         (\bar{Q}_i\cdot\bar{Q}_i)_{t}=-(\bar{Q}_i\cdot\operatorname{div}(\omega))=-(2g-2)(\bar{Q}_i\cdot\bar{\infty})-(\bar{Q}_i\cdot V).   
 \end{equation*} 
In order to compute the intersection multiplicity $(\bar{Q}_i\cdot V)$, we assume that $(\bar{\infty},V)=0$. Otherwise we can add a suitable multiple of the special fiber to $V$. Since $\bar{\infty}$ meets $\Gamma_\infty$ transversely and meets no other components in the special fiber, we have $\Gamma_\infty\varsubsetneq|V|$. 

A canonical $\mathbb{Q}$-divisor $\mathcal{K}$ is a $\mathbb{Q}$-divisor satisfying $\mathcal{O}_\mathcal{C}(\mathcal{K})\cong\boldsymbol{\omega}_\mathcal{C}$, where $\boldsymbol{\omega}_\mathcal{C}$ is the canonical sheaf of $\mathcal{C}$. We now show that $\operatorname{div}(\omega)$ is a canonical $\mathbb{Q}$-divisor $\mathcal{K}$. Since $\operatorname{div}(\frac{dx}{2y})=(2g-2)\infty$ is a canonical divisor on $C$,\cite[Proposition~2.5]{CurillaArithmetic} implies that there is a canonical $\mathbb{Q}$-divisor $\mathcal{K}$ on $\mathcal{C}$ that extends $\operatorname{div}(\frac{dx}{2y})$. The divisor $\mathcal{K}$ is unique up to a multiple of the special fiber $\mathcal{C}_s$. In order to fix $\mathcal{K}$, we choose $\mathcal{K}$ satisfying the condition $\Gamma_\infty\varsubsetneq|\mathcal{K}-(2g-2)\bar{\infty}|$. Since the canonical sheaf $\boldsymbol{\omega}_\mathcal{C}$ is not isomorphic to the relative cotangent bundle $\Omega_{\mathcal{C}}$ in general, the divisor $\operatorname{div}(\omega)$ might not be equal to $\mathcal{K}$. Since the reduction is semistable, $\boldsymbol{\omega}_{\mathcal{C}}$ is a line bundle by \cite[\S1]{Deligne}. It is shown in \cite[Corollary~6.4.13]{Liualgebraic} that the isomorphism $\boldsymbol{\omega}_\mathcal{C}\cong\Omega_{\mathcal{C}}$ holds on the smooth locus. This means that the divisor $\operatorname{div}(\omega)-\mathcal{K}$ is only supported on the vertical divisors of multiplicity at least $2$ by \cite[\S5.1]{BalakrishnanQuadratic}. Then, we have $\mathcal{K}=(2g-2)\bar{\infty}+V$ because the multiplicity of any irreducible component of $\mathcal{C}_s$ is $1$.
 
For a $\mathbb{Q}$-divisor $E=\sum e_\Gamma\Gamma\in\mathbb{Q}\operatorname{Div}(\mathcal{C}_s)$, the divisor $\sum d_\Gamma\Gamma$ is denoted by $E_1$ where $d_\Gamma=(E\cdot\Gamma)$. We obtain $e_\Gamma=-L^+(E_1,\Gamma)$ by \cite[\S6]{Bianchicoleman}. Writing $V=\sum a_\Gamma\Gamma$, we have the following formula for $(\bar{Q}_i\cdot V)$ in terms of $L^+$ and intersection multiplicities \[a_{\Gamma_{Q_i}}=(\bar{Q}_i\cdot V)=-L^+(\Gamma_{Q_i},K_G-(2g-2)\Gamma_\infty)=-L^+(\Gamma_{Q_i},K_G)+(2g-2)L^+(\Gamma_{Q_i},\Gamma_\infty).\]
We have $ (\bar{Q}_i\cdot\bar{Q}_i)_{t}=-(2g-2)(\bar{Q}_i\cdot\bar{\infty})+L^+(K_G,\Gamma_{Q_i})-(2g-2)L^+(\Gamma_{Q_i},\Gamma_\infty)$. We obtain the following formula
 \begin{equation*}
     \begin{aligned}
          &\lambda_N(a)\\
          \equiv\quad&(2-4g)\sum_{i=1}^g(\bar{Q}_i\cdot\bar{\infty})+2\sum_{1\leq i<j\leq g}(\bar{Q_i}\cdot\bar{Q_j})_t+\sum_{i=1}^gL^+(K_G,\Gamma_{Q_i})+(2-2g)\sum_{i=1}^gL^+(\Gamma_{Q_i},\Gamma_\infty)\\
          &+L^+(\sum_{i=1}^g\Gamma_{Q_i}-g\Gamma_\infty,\sum_{i=1}^g\Gamma_{Q_i}-g\Gamma_\infty)\\
         \equiv\quad&(2-4g)(\sum_{i=1}^g\bar{Q}_i\cdot\bar{\infty})+2\sum_{1\leq i<j\leq g}(\bar{Q_i}\cdot\bar{Q_j})_t+\sum_{i=1}^gL^+(K_G,\Gamma_{Q_i})+(2-4g)L^+(\sum_{i=1}\Gamma_{Q_i},\Gamma_\infty)\\
          &+\sum_{i=1}^gL^+(\Gamma_{Q_i},\Gamma_{Q_i})+2\sum_{1\leq i<j\leq g}L^+(\Gamma_{Q_i},\Gamma_{Q_j}).
     \end{aligned}
 \end{equation*}
The proof is complete.
\end{proof}
We obtain different formulas for $\lambda_{2\Theta}(a)$ according to the parity of the genus. Firstly, we consider the \textit{odd} case, the \textit{even} case is proved by a similar method.      
\begin{lemma}
    Suppose now that $g=2n-1$. As in \S\ref{ch6}, let $D_a\coloneqq Q_1+\cdots+Q_g-g\infty\in\operatorname{Div}^0(C)$ be the representation for $a\in J(C)(K)\setminus|\Theta|$ where $Q_i\in C(\overline{K})\setminus\{\infty\}$ and $\iota(Q_i)\ne Q_j$ for all $i\ne j$. Then, we have the following formula 
    \begin{equation}\label{odd}
    \begin{aligned}
        &\lambda_{2\Theta}(a)\\
        \equiv\quad&2\sum_{i=1}^{2n-1}(\bar{Q_i}\cdot\iota\bar{Q}_{i+1}+\cdots+\iota\bar{Q}_{i+n-1}+\bar{\infty}
	 )-\frac{n}{2n-1}\sum_{i=1}^{2n-1}L^+(\Gamma_{Q_i},\Gamma_{Q_i})-\frac{n-1}{2n-1}\sum_{i=1}^{2n-1}L^+(\Gamma_{\iota Q_i},\Gamma_{\iota Q_i})\\
     &-\frac{1}{2n-1}\sum_{i=1}^{2n-1}L^+(K_G,n\Gamma_{Q_i}+\Gamma_{\iota Q_{i+1}}+\cdots+\Gamma_{\iota Q_{i+n-1}}+\Gamma_\infty)\\
     &+2\sum_{i=1}^{2n-1}L^+(\Gamma_{Q_i},\Gamma_{\iota Q_{i+1}}+\cdots+\Gamma_{\iota Q_{i+n-1}}+\Gamma_\infty).\\
    \end{aligned} 
    \end{equation}
\end{lemma}
\begin{proof}
Let $a=[Q_1+...+Q_{2n-1}-(2n-1)\infty]\in J(C)(K)\setminus|\Theta|$. Theorem~\ref{Thm6.2} implies that   
\begin{equation}\label{odd1}
    \lambda_{2\Theta}(a)=\sum_{1\leq i< j\leq 2n-1}2\langle Q_i,\iota Q_j\rangle_{\operatorname{ad}}+2\sum_{i=1}^{2n-1}\langle Q_i,\infty\rangle_{\operatorname{ad}}.
\end{equation}
By rearranging the order of $\{Q_i\}$ in the pairings $\{
	\langle Q_i,\iota Q_j\rangle_{\operatorname{ad}}\}$, we rewrite $\lambda_{2\Theta}(a)$ as follows 
	 \[\lambda_{2\Theta}(a)=2\sum_{i=1}^{2n-1}\langle Q_i,\iota Q_{i+1}
	 +\cdots+\iota Q_{i+n-1}+\infty\rangle_{\operatorname{ad}}\] 
where $\iota Q_i=\iota Q_{i\operatorname{mod} 2n-1}$ when $i>2n-1$. By the definition of admissible pairings, we have
\begin{equation*}
     \begin{aligned}
         &\quad\lambda_{2\Theta}(a)\\
         &\equiv2\sum_{i=1}^{2n-1}(\bar{Q_i}\cdot\iota\bar{Q}_{i+1}+..+\iota\bar{Q}_{i+n-1}+\bar{\infty})+2\sum_{i=1}^{2n-1}G(\Gamma_{Q_i},\Gamma_{\iota Q_{i+1}}+...+\Gamma_{\iota Q_{i+n-1}}+\Gamma_\infty)\\
	 &\equiv2\sum_{i=1}^{2n-1}(\bar{Q_i}\cdot\iota\bar{Q}_{i+1}+..+\iota\bar{Q}_{i+n-1}+\bar{\infty}
	 )-\frac{n}{2n-1}\sum_{i=1}^{2n-1}L^+(\Gamma_{Q_i},\Gamma_{Q_i})-\frac{n-1}{2n-1}\sum_{i=1}^{2n-1}L^+(\Gamma_{\iota Q_i},\Gamma_{\iota Q_i})\\
     &\quad-\frac{1}{2n-1}\sum_{i=1}^{2n-1}L^+(K_G,n\Gamma_{Q_i}+\Gamma_{\iota Q_{i+1}}+...+\Gamma_{\iota Q_{i+n-1}}+\Gamma_\infty)\\
     &\quad+2\sum_{i=1}^{2n-1}L^+(\Gamma_{Q_i},\Gamma_{\iota Q_{i+1}}+..+\Gamma_{\iota Q_{i+n-1}}+\Gamma_\infty)\\
     \end{aligned}  
\end{equation*} 
The proof is complete. 
\end{proof}
\begin{lemma}
We assume that $g=2n$ is even. For any $a\in J(C)(K)\setminus|\Theta|$, the representation $a=[D_a]$ is defined in \S\ref{ch6}. Then, the function $\lambda_{2\Theta}(a)$ is given by 
  \begin{equation}\label{even}
      \begin{aligned}
          \lambda_{2\Theta}(a)&\equiv2\sum_{i=1}^{2n}(\bar{Q_i}\cdot\iota\bar{Q}_{i+1}+..+\iota\bar{Q}_{i+n-1}+\bar{\infty}
	 )+2\sum_{i=1}^{n}(\bar{Q}_i\cdot\iota\bar{Q}_{i+n})\\
     &\quad-\frac{1}{2n}\sum_{i=1}^{2n}L^+(K_G,n\Gamma_{Q_i}+\Gamma_{\iota Q_{i+1}}+...+\Gamma_{\iota Q_{i+n-1}}+\Gamma_\infty)-\frac{1}{2}\sum_{i=1}^{2n}L^+(\Gamma_{Q_i},\Gamma_{Q_i})\\
     &\quad-\frac{n-1}{2n}\sum_{i=1}^{2n}L^+(\Gamma_{\iota Q_i},\Gamma_{\iota Q_i})+2\sum_{i=1}^{2n}L^+(\Gamma_{Q_i},\Gamma_{\iota Q_{i+1}}+..+\Gamma_{\iota Q_{i+n-1}}+\Gamma_\infty)\\
     &\quad-\frac{1}{2n}\sum_{i=1}^{n}L^+(K_G,\Gamma_{Q_i}+\Gamma_{\iota Q_{i+n}})-\frac{1}{2n}\sum_{i=1}^{n}L^+(\Gamma_{Q_i},\Gamma_{Q_i})-\frac{1}{2n}\sum_{i=1}^{n}L^+(\Gamma_{\iota Q_{i+n}},\Gamma_{\iota Q_{i+n}})\\
     &\quad+2\sum_{i=1}^{n}L^+(\Gamma_{Q_i},\Gamma_{\iota Q_{i+n}})
      \end{aligned}
  \end{equation}
\end{lemma}
\begin{proof}
As in the odd genus case, we also need to rearrange the order of $\{Q_i\}$ in Theorem~\ref{Thm6.2}. It can be rewritten as follows   
\[\lambda_{2\Theta}(a)=2\sum_{i=1}^{2n}\langle Q_i,\iota Q_{i+1}+...+\iota Q_{i+n-1}\rangle_{\operatorname{ad}}+2\sum_{i=1}^{n}\langle Q_i,\iota Q_{i+n}\rangle_{\operatorname{ad}}+2\sum_{i=1}^{2n}\langle Q_i,\infty\rangle_{\operatorname{ad}}\]
where $Q_i=Q_{(i)\operatorname{mod} 2n}$ when $i>2n$. Note that there are extra terms $2\sum_{i=1}^{n}\langle Q_i,\iota Q_{i+n}\rangle_{\operatorname{ad}}$ in the above formula compared with the odd genus case \eqref{odd1}. The following formula for $\lambda_{2\Theta}(a)$ is obtained similarly to the odd genus case,
\begin{align*}
\lambda_{2\Theta}(a)
&\equiv 2\sum_{i=1}^{2n}(\bar{Q_i}\cdot\iota\bar{Q}_{i+1}+\cdots+\iota\bar{Q}_{i+n-1})+2\sum_{i=1}^{2n}G(\Gamma_{Q_i},\Gamma_{\iota Q_{i+1}}+\cdots+\Gamma_{\iota Q_{i+n-1}})\\
	   &\quad+2\sum_{i=1}^{n}(\bar{Q}_i\cdot\iota\bar{Q}_{i+n})+2\sum_{i=1}^{n}G(\Gamma_{Q_i},\Gamma_{\iota Q_{i+n}})+2\sum_{i=1}^{2n}(\bar{Q_i}\cdot\bar{\infty})+2\sum_{i=1}^{2n}G( \Gamma_{Q_i},\Gamma_\infty)\\
     &\equiv 2\sum_{i=1}^{2n}(\bar{Q_i}\cdot\iota\bar{Q}_{i+1}+\cdots+\iota\bar{Q}_{i+n-1}+\bar{\infty}
	 )+2\sum_{i=1}^{n}(\bar{Q_i}\cdot\iota\bar{Q}_{i+n})\\
     &\quad-\frac{1}{2n}\sum_{i=1}^{2n}L^+(K_G,n\Gamma_{Q_i}+\Gamma_{\iota Q_{i+1}}+\cdots+\Gamma_{\iota Q_{i+n-1}}+\Gamma_\infty)-\frac{1}{2}\sum_{i=1}^{2n}L^+(\Gamma_{Q_i},\Gamma_{Q_i})\\
     &\quad-\frac{n-1}{2n}\sum_{i=1}^{2n}L^+(\Gamma_{\iota Q_i},\Gamma_{\iota Q_i})+2\sum_{i=1}^{2n}L^+(\Gamma_{Q_i},\Gamma_{\iota Q_{i+1}}+\cdots+\Gamma_{\iota Q_{i+n-1}}+\Gamma_\infty)\\
     &\quad-\frac{1}{2n}\sum_{i=1}^{n}L^+(K_G,\Gamma_{Q_i}+\Gamma_{\iota Q_{i+n}})-\frac{1}{2n}\sum_{i=1}^{n}L^+(\Gamma_{Q_i},\Gamma_{Q_i})-\frac{1}{2n}\sum_{i=1}^{n}L^+(\Gamma_{\iota Q_{i+n}},\Gamma_{\iota Q_{i+n}})\\
     &\quad+2\sum_{i=1}^{n}L^+(\Gamma_{Q_i},\Gamma_{\iota Q_{i+n}}).\\
     \end{align*}
This completes the proof.
\end{proof}
We now prove Theorem~\ref{main}. 
\begin{proof}
By the definition of $\lambda_N$, we obtain the following formula for $-\sum_{1\leq i<j\leq g}\log|f_j(Q_i-\infty,t)|$
\begin{equation*}
    \begin{aligned}
    &-\sum_{1\leq i<j\leq g}\log|f_j(Q_i-\infty,t)|\\
    &=\sum_{1\leq i<j\leq g}\langle Q_i-\infty,Q_j+\iota Q_j-2\infty\rangle_{N,t}\\
          &=\sum_{1\leq i<j\leq g}(\bar{Q_i}-\bar{\infty}\cdot\bar{Q_j}+\bar{\iota Q_j}-2\bar{\infty})_t-(\Phi(Q_i-\infty)\cdot\Phi(Q_j+\iota Q_j-2\infty))\\
          &=\sum_{1\leq i<j\leq g}(\bar{Q_i}-\bar{\infty}\cdot\bar{Q_j}+\bar{\iota Q_j}-2\bar{\infty})_t+L^+(\Gamma_{Q_i}-\Gamma_\infty,\Gamma_{Q_j}+\Gamma_{\iota Q_j}-2\Gamma_\infty).\\
    \end{aligned}
\end{equation*}
We have $f_j(Q_i-\infty,t)^2=(x(Q_i)-x(Q_j))^2$ since the divisor $\operatorname{div}(x-x(Q_j))=Q_j+\iota Q_j-2\infty$ is principal. Therefore, there is an equality 
\begin{equation}\label{fj}
    \log|f_j(Q_i-\infty,t)| =\log|f_i(Q_j-\infty,t)|
\end{equation} for any $i\neq j$. 

Suppose that $g=2n-1$ is \textit{odd}. We rearrange the distribution of $f_j(Q_i-\infty,t)^2$ in the sum $\sum_{1\leq i<j\leq g}\log|f_j(Q_i-\infty,t)|$ by \eqref{fj}. We deduce
    \begin{equation}
        \begin{aligned}\label{fiodd}
        &\quad-2\sum_{1\leq i<j\leq g}\log|f_j(Q_i-\infty,t)|\\
        &=2\sum_{i=1}^{2n-1}\langle Q_i-\infty,\sum_{j=i+1}^{i+n-1}Q_{j}+\iota Q_{j}-2(n-1)\infty\rangle_{N,t}\\
       &\equiv\sum_{i=1}^{2n-1}\left(2(\bar{Q}_i\cdot\iota\bar{Q}_{i+1}+...+\iota \bar{Q}_{i+n-1}+(3-3n)\bar{\infty})-2(\iota\bar{Q}_i\cdot(n-1)\bar{\infty})\right)+2\sum_{1\leq i<j\leq g}(\bar{Q_i}\cdot\bar{Q_j})_t\\
        &\quad+\sum_{i=1}^{2n-1}\left(2L^+(\Gamma_{Q_i},\Gamma_{\iota Q_{i+1}}+...+\Gamma_{\iota Q_{i+n-1}}+(3-3n)\Gamma_\infty)+2L^+(\Gamma_{\iota Q_i},(1-n)\Gamma_{\infty}\right)\\
        &\quad+2\sum_{1\leq i<j\leq g}L^+(\Gamma_{Q_i},\Gamma_{Q_j}).\\
        \end{aligned}
    \end{equation}
Combining \eqref{lambda},\eqref{odd}, and \eqref{fiodd}, we obtain the following formula 
\begin{equation*}
    \begin{aligned}
    &\quad\lambda_N(a)+\lambda_{2\Theta}(a)+2\sum_{1\leq i<j\leq g}\log|f_j(Q_i-\infty,t)|\\
        &\equiv(1-n)\sum_{i=1}^{2n-1}\bigl(2(\bar{Q}_i-\iota\bar{Q}_i\cdot\bar{\infty})+2L^+(\Gamma_{Q_i}-\Gamma_{\iota Q_i},\Gamma_\infty)-\frac{1}{2n-1}L^+(\Gamma_{Q_i},\Gamma_{Q_i})\\
        &\quad+\frac{1}{2n-1}L^+(\Gamma_{\iota Q_i},\Gamma_{\iota Q_i})-\frac{1}{2n-1}L^+(K_G,\Gamma_{ Q_i}-\Gamma_{\iota Q_{i}})\bigr).\\
        \end{aligned}
\end{equation*}
It is easy to compute that the sum of all $L^+$ pairings in the above equality is equal to $2(1-n)\sum_{i=1}^{2n-1}(\iota \bar{Q}_i-\bar{Q}_i,\bar{\infty})$  up to a constant by Lemma~\ref{lemma difference}. Hence, the right hand side of the equality is a constant which proves the \textit{odd} genus case of Theorem~\ref{main}.

We now consider the \textit{even} case. Following the proof of the \textit{odd} case, we also need to swap some terms in $2\sum_{1\leq i<j\leq g}\log|f_j(Q_i-\infty,t)|$. By doing so, we obtain 
\begin{equation}\label{fieven}
    \begin{aligned}
    &\quad-2\sum_{1\leq i<j\leq g}\log|f_j(Q_i-\infty,t)|\\
    &=2\sum_{i=1}^{2n}\langle Q_i-\infty,\sum_{j=i+1}^{i+n-1}Q_{j}+\iota Q_{j}-2(n-1)\infty\rangle_{N,t}+2\sum_{i=1}^n\langle Q_i-\infty,Q_{i+n}+\iota Q_{i+n}-2\infty\rangle_{N,t}\\
    &\equiv\sum_{i=1}^{2n}2\left((\bar{Q}_i\cdot\iota \bar{Q}_{i+1}+...+\iota \bar{Q}_{i+n-1}+(3-3n)\bar{\infty})-(\iota\bar{Q}_i\cdot(n-1)\bar{\infty})\right)+2\sum_{1\leq i<j\leq g}(\bar{Q_i}\cdot\bar{Q_j})_t\\
    &\quad+2\sum_{i=1}^{2n}\left(L^+(\Gamma_{Q_i},\Gamma_{\iota Q_{i+1}}+...+\Gamma_{\iota Q_{i+n-1}}+(3-3n)\Gamma_\infty)+L^+(\Gamma_{\iota Q_i},(1-n)\Gamma_{\infty}\right)\\
    &\quad+2\sum_{1\leq i<j\leq g}L^+(\Gamma_{Q_i},\Gamma_{Q_j})+2\sum_{i=1}^n(\bar{Q}_i\cdot\iota \bar{Q}_{i+n}-2\infty)-2\sum_{i=1}^n(\iota\bar{Q}_{i+n}+\bar{Q}_{i+n}\cdot\infty)\\
    &\quad+2\sum_{i=1}^nL^+(\Gamma_{Q_i},\Gamma_{\iota Q_{i+n}}-2\Gamma_\infty)-2\sum_{i=1}^nL^+(\Gamma_{\iota Q_{i+n}}+\Gamma_{Q_{i+n}},\Gamma_\infty).
    \end{aligned}
\end{equation}
The following formula for $\lambda_N(a)+\lambda_{2\Theta}(a)+2\sum_{1\leq i<j\leq g}\log|f_j(Q_i-\infty,t)|$ is obtained by combining \eqref{lambda},\eqref{even} and \eqref{fieven}  
\begin{equation*}
    \begin{aligned}
    &\quad\lambda_N(a)+\lambda_{2\Theta}(a)+2\sum_{1\leq i<j\leq g}\log|f_j(Q_i-\infty,t)|\\
    &\equiv\sum_{i=1}^n2(\iota\bar{Q}_i-\bar{Q}_i\cdot(n-1)\bar{\infty})+2(\iota\bar{Q}_{i+n}-\bar{Q}_{i+n}\cdot n\bar{\infty})\\
        &\quad+\sum_{i=1}^n2L^+(\Gamma_{\iota Q_{i+n}}-\Gamma_{Q_{i+n}},n\Gamma_\infty)+2L^+(\Gamma_{\iota Q_i}-\Gamma_{Q_i},(n-1)\Gamma_\infty)\\
        &\quad+\frac{(n-1)}{2n}\sum_{i=1}^{n}L^+(K_G,\Gamma_{Q_i}-\Gamma_{\iota Q_i})+\frac{1}{2}\sum_{i=1}^{n}L^+(K_G,\Gamma_{Q_{i+n}}-\Gamma_{\iota Q_{i+n}})\\
        &\quad+\frac{1}{2}\sum_{i=1}^{n}L^+(\Gamma_{Q_{i+n}},\Gamma_{Q_{i+n}})-L^+(\Gamma_{\iota Q_{i+n}},\Gamma_{\iota Q_{i+n}})+\frac{n-1}{2n}\sum_{i=1}^{n}(L^+(\Gamma_{Q_i},\Gamma_{Q_i})-L^+(\Gamma_{\iota Q_i},\Gamma_{\iota Q_i})).\\
    \end{aligned}
\end{equation*}
It is also true that the sum of all $L^+$ pairings cancels out the sum of all intersection multiplicities in the above formula by Lemma~\ref{lemma difference}. The \textit{even} genus case is proved. 
\end{proof}
\begin{rem}\label{conjecture}
We expect that Theorem~\ref{main} still holds over archimedean local fields. When $C$ is a curve over $\mathbb{Q}$ or an imaginary quadratic field, it is true because of the following.

Let $J(C)$ be the Jacobian of a curve $C$ over a number field $F$. Let $D$ and $E$ be divisors of degree zero on $C$, and let $h_{2\Theta}(\,,\,)$ be the bilinear pairing induced from the canonical height associated to $2\Theta$. The following equality is due to Faltings and Hriljac, see \cite[Theorem~4]{FaltingsCalculus} and \cite[Theorem~3.1]{Hriljacheights}. 
\[\sum_v\langle D,E\rangle_{v,N}=-h_{2\Theta}([D],[E])\] where the sum is over all places of $F$. When $D=E$, by fixing a non-zero set-theoretic section $t$ of the tangent bundle, defined over $F$. The equality is rewritten as \[\sum_v\langle D,D\rangle_{v,N}=-h_{2\Theta}([D])\] see \cite[\S5]{Grosslocalheights}. The sum $\sum_v\langle D,D\rangle_{v,N}$ is independent of the choice of $t$. Moreover,\cite[Theorem~9.3]{SilvermanDiophantine} shows that the canonical height admits the following decomposition into a sum of N\'eron functions \[\sum_vn_v\lambda_{2\Theta}([D])=h_{2\Theta}([D])\] where $n_v=[F_v:\mathbb{Q}_v]/[F:\mathbb{Q}]$.

Now let $C$ denote an odd-degree hyperelliptic curve of genus $g$ over a number field $F$, and $a=[D_a]\in J(C)(K)\setminus|\Theta|$. Theorem \ref{main} shows that \eqref{maineq} holds over non-archimedean places. Combining the above discussion with the product formula, we have the equality \[\sum\lambda_{N,v}(a)=\sum\left(-\lambda_{2\Theta}(a)-2\sum_{1\leq i<j\leq g}\log|f_j(Q_i-\infty,t)|_v+c_v\right)\] where the sum is over all archimedean places. When $F$ is a number field with only one archimedean place, namely $F=\mathbb{Q}$ or $F$ is an imaginary quadratic field, Theorem~\ref{main} follows for the archimedean place.              
\end{rem}
\subsection{Coleman--Gross heights and $p$-adic N\'eron functions}\label{ch2.4}
Let $C$ be a curve over a number field $F$. For a finite place $v$ of $F$ and a prime number $p$, let $\chi_v$ be a continuous homomorphism $F_v^\times\to\mathbb{Q}_p$, and let $\pi_v$ be a uniformizer of the ring of integers $\mathcal{O}_v$. The character $\chi_v$ is called \textit{ramified} if $\chi_v(\mathcal{O}_v^\times)\ne 0$, where $\mathcal{O}_v^\times$ is the unit group of $\mathcal{O}_v$, in which case necessarily $v\mid p$. Otherwise $\chi_v$ is called \textit{unramified}.

The \textit{Coleman--Gross height} $h^{\operatorname{CG}}$ is a bilinear pairing $\operatorname{Div}^0(C)\times\operatorname{Div}^0(C)\to\mathbb{Q}_p$ constructed in \cite{ColemanGrossp-adic}. It can be written as the sum of local height pairings $h^{\operatorname{CG}}_{\chi_v}$ with respect to an id\'ele class character $\chi$. Here, $\chi$ is a collection of continuous homomorphisms $\chi_v$, one for each finite place $v$ of $F$, satisfying an analogue of the product formula. The local Coleman--Gross height pairing $h^{\operatorname{CG}}_{\chi_v}\colon\operatorname{Div}^0(C_{F_v})\times\operatorname{Div}^0(C_{F_v})\to\mathbb{Q}_p$ with respect to $\chi_v$ is defined in \cite{ColemanGrossp-adic} for divisors with disjoint support and extended to all divisors in \cite[\S2]{BalakrishnanBesser}, depending on a choice of nonzero $F_v$-section $t$ of the tangent bundle as in \S\ref{ch2.2}. Suppose that $\chi_v$ is unramified, then for any divisors $D,E\in\operatorname{Div}^0(C_{F_v})$, we have \[h^{\operatorname{CG}}_{\chi_v,t}(D,E)=\chi_v(\pi_v)\cdot\alpha_t(D,E),\]where $\alpha_t(\,,\,)$ is the same pairing in Proposition~\ref{extendedpairing}. Although the N\'eron pairing and the local Coleman--Gross height pairing depend on the choice of $t$, their sums over all places (i.e. the canonical N\'eron height pairing and the Coleman--Gross height pairing) are independent of this choice.   

Moreover, for any $a\in J(C)(F_v)\setminus|\Theta|$, the \textit{local Coleman--Gross height function} $\lambda_{\chi_v}^{\operatorname{CG}}(a)$ is defined as $h^{\operatorname{CG}}_{\chi_v}(D_a,D_a)$. 


Let $A$ be an abelian variety over $F$, and let $D\in\operatorname{Div}(A)$. The \textit{$p$-adic N\'eron function} $\lambda_{D,\chi_v}$ with respect to $D$ and $\chi_v$ is introduced in \cite[\S4]{Bianchicoleman} with values in $\mathbb{Q}_p$. When $\chi_v$ is unramified, it is shown in \cite[Remark~4.12]{Bianchicoleman} that there is a function $\beta\colon A(F_v)\setminus|D|\to\mathbb{Q}$ such that the $p$-adic N\'eron function can be written as \[\lambda_{D,\chi_v}(a)=-\chi_v(\pi_v)\cdot\beta(a),\] and the N\'eron function $\lambda_D$ can be expressed as \[\lambda_D(a)=\log q_v\cdot\beta(a).\] 

So our results for N\'eron functions and N\'eron pairings immediately imply the corresponding results 
for their $p$-adic analogues. It is shown in \cite[Theorem~6.3]{Bianchicoleman} that unramified $p$-adic N\'eron functions agree with unramified local Coleman--Gross height functions on curves of genus $2$. The following theorem generalizes this result and follows from Theorem~\ref{main}.
\begin{theo}\label{padiccomparision}
Let $C$ be an odd-degree hyperelliptic curve of genus $g$ over a non-archimedean local field $K$. Let $\lambda_{2\Theta,\chi_v}$ be the $p$-adic N\'eron function with respect to an unramified character $\chi_v$. There is the following identity up to an additive constant \[\lambda_{\chi_v}^{\operatorname{CG}}(a)\equiv
\lambda_{2\Theta,\chi_v}(a)-2\sum_{1\leq i<j\leq g}\chi_v(f_j(Q_i-\infty,t))\]
where $a=[D_a]\in J(C)(K)\setminus|\Theta|$.
\end{theo}
 
\section{Explicit formulas and algorithms for $T(a)$}\label{ch8}
Let $K$ be a non-archimedean local field, and let $R$ be its valuation ring with residue field $k$. Let $C$ be an odd-degree hyperelliptic curve of genus $g$ over $K$. Theorem~\ref{Thm6.2} shows that the N\'eron function associated to the divisor $\Theta$ is given by \[\lambda_{\Theta}(a)=\sum_{1\leq i<j\leq g}\langle Q_i,\iota Q_j\rangle_{\operatorname{ad}}+\sum_{i=1}^g\langle Q_i,\infty\rangle_{\operatorname{ad}}\] for any point $a=[D_a]\in J(C)(K)\setminus|\Theta|$.  
\begin{defi}\label{T}
We define a function $I\colon J(C)(K)\setminus|\Theta|\to\mathbb{R}$ by \[I(a)=\sum_{1\leq i<j\leq g}(Q_i\cdot\iota Q_j)+\sum_{i=1}^g(Q_i\cdot\infty)\] where $(\,\cdot\,)$ is the intersection multiplicity defined in \S\ref{3.2}. We also define a function $T\colon J(C)(K)\setminus|\Theta|\to\mathbb{R}$ by \[T(a)=\sum_{1\leq i<j\leq g}G(\operatorname{sp}(Q_i),\operatorname{sp}(\iota Q_j))+\sum_{i=1}^gG(\operatorname{sp}(Q_i),\operatorname{sp}(\infty))\] where $G(\,,\,)$ is the Green function defined in \S\ref{3.2}.  
\end{defi}
\begin{coro}\label{decomposition}
    There is a decomposition of the N\'eron function as follows \[\lambda_{\Theta}(a)=I(a)+T(a).\] 
\end{coro}  

\begin{proof}
Recall the formula $\langle P,Q\rangle_{\operatorname{ad}}=i(P,Q)+G(\operatorname{sp}(P),\operatorname{sp}(Q))$ from \S\ref{3.2}. Combining it with Theorem~\ref{Thm6.2}, we have the result. 
\end{proof}

\begin{rem}
Suppose $C$ is a curve of genus $2$ over $K$ and $a=[P_1-P_2]\in J(C)(K)\setminus|\Theta|$. Corollary~\ref{Bianchi} shows that $T(a)$ can also be defined by $G(\operatorname{sp}(P_1),\operatorname{sp}(P_2))+G(\operatorname{sp}(P_1)+\operatorname{sp}(P_2),\operatorname{sp}(\infty))$. It is equivalent to the formula for $T$ in Definition~\ref{T}.
\end{rem}

\begin{rem}
Suppose that $C$ has semistable reduction over $K$. Let $L/K$ be a finite ramified extension of ramification index $e>1$, hence $C$ also has semistable reduction over $L$. Let $P\in C(L)$, and let $G_{C_K}$ be the dual graph of the minimal regular model $\mathcal{C}$ of $C$ over $K$. The closure of $P$ may meet the special fiber of $\mathcal{C}_{\mathcal{O}_L}$ at a singular point. After base change to $\mathcal{O}_L$, we obtain a minimal regular model $\mathcal{C}_L$ of $C_L$ by desingularization. The dual graph $G_{C_L}$ of $\mathcal{C}_L$ is obtained by subdividing every edge of $G_{C_K}$ into $e$ edges. Then the point $P$ specializes to a vertex of $G_{C_L}$. 
\end{rem}

A point $a=[D_a]\in J(C)(K)\setminus|\Theta|$ is called \textit{split unramified} if all $Q_i$ lie in $C(K^{\operatorname{ur}})$, where $K^{\operatorname{ur}}$ is the maximal unramified extension of $K$. Then the specialization of each $Q_i$ lies in the set of vertices of $G_C$. In this section, we determine the possible values of $T$ at split unramified points when the minimal proper regular model of $C$ has certain semistable reduction types. To this end, we use the formula in Lemma~\ref{kayagreen}, which expresses the Green function on vertices. We first need to compute the Green function. 
\subsection{Explicit formula for the Green function}
A complete classification of the special fibers of the proper regular minimal models of genus $2$ curves was given by Namikawa and Ueno in \cite{Uenopencil}. In particular, there are seven different configurations of the special fibers for genus $2$ curves with stable model. They are denoted by $\mathrm{I}$, $\mathrm{II}(n)$, $\mathrm{III}(n)$, $\mathrm{IV}(n,m)$, $\mathrm{V}(n,m)$, $\mathrm{VI}(n,m,l)$, $\mathrm{VII}(n,m,l)$ in \cite{MoriwakiBogomolov} where the positive integers $n,m,l$ denote the respective thickness of singular points.
     
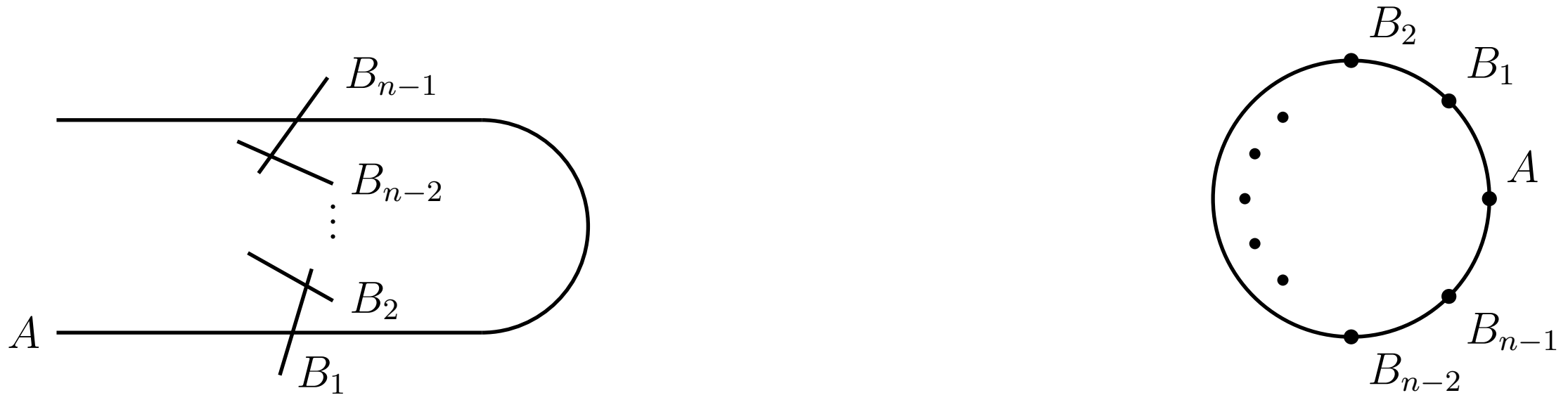
\begin{figure}[htbp]
 \centering
\begin{minipage}[t]{0.45\textwidth}   
   \raggedright
     \begin{tikzpicture}[line width=1pt]
\draw (0,2) -- (4,2);
\draw (4,2) arc[start angle=90,end angle=-90,radius=1];
\draw (4,0) -- node[pos=1,left] {$A$}(0,0);
\draw
(2.55,2.4)  
-- node[pos=0,right] {$B_{n-1}$}(1.9,1.5);
\draw
(1.7,1.8)    
--node[pos=1,right] {$B_{n-2}$} (2.6,1.4);
\node (dots) at (2.6,1.15) {$\vdots$};
\draw
(1.8,0.75)
--node[pos=1,right] {$B_2$}(2.60,0.3);
\draw
(2.4,0.6)
-- node[pos=1,right] {$B_1$}(2.1,-0.4);   
\end{tikzpicture}
\end{minipage}
\hfill
\begin{minipage}[t]{0.45\textwidth}
\centering
\begin{tikzpicture}[line width=1pt]
\draw (0,0) circle (1.3);
\fill (0:1.3)  circle (2pt);
    \node[above right] at (0:1.3) {$A$};  
\fill (45:1.3)  circle (2pt);
    \node[above right] at (45:1.3) {$B_1$};
\fill (90:1.3)  circle (2pt);
    \node[above right] at (90:1.3) {$B_2$};
\fill (130:1)  circle (1.5pt);
\fill (155:1)  circle (1.5pt);
\fill (180:1)  circle (1.5pt);
\fill (205:1)  circle (1.5pt);
\fill (230:1)  circle (1.5pt);
\fill (270:1.3)  circle (2pt);
    \node[below right] at (270:1.3) {$B_{n-2}$};
\fill (315:1.3)  circle (2pt);
    \node[below right] at (315:1.3) {$B_{n-1}$};    
\end{tikzpicture}   
\end{minipage}
\caption{The special fiber of reduction type $\mathrm{III}(n)$ and its dual graph}
    \label{fig:III}
\end{figure}

Let $C$ be a curve of genus $2$ with reduction type $\mathrm{III}(n)$. The special fiber of its stable model is an elliptic curve with a node of thickness $n$. Desingularizing the node by a chain of $n$ rational curves, the special fiber of the minimal regular model consists of the elliptic curve whose node is replaced by the chain. More generally, we say that a curve $C$ of genus $g$ has reduction type $\mathrm{III}(n)$ if it has semistable reduction and if the special fiber of its stable model is an irreducible curve of genus $g-1$ which has a node of thickness $n$. The special fiber of the minimal regular model is obtained by replacing the node by a chain of $n-1$ rational components, see Figure\ref{fig:III}.

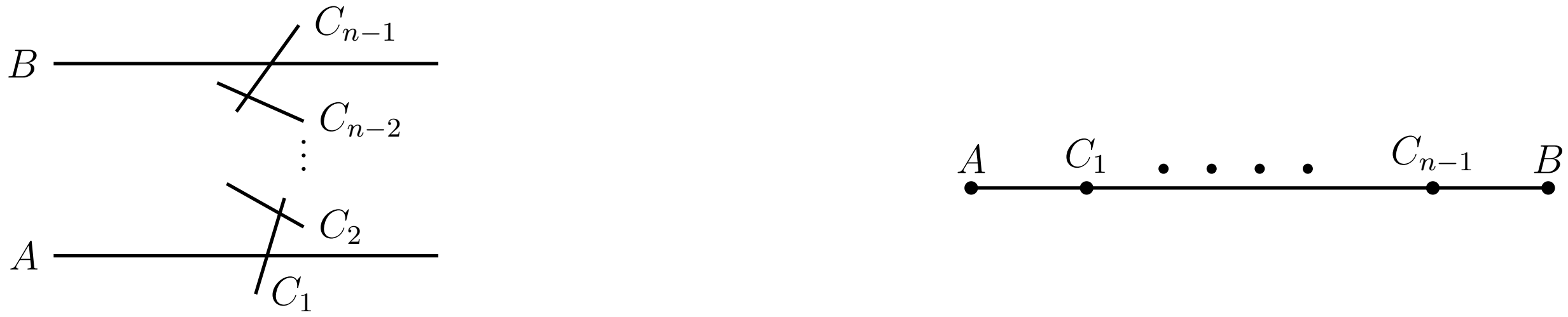
\begin{figure}[htbp]
 \centering
\begin{minipage}[c]{0.45\textwidth}   
   \raggedright
     \begin{tikzpicture}[line width=1pt]
\draw (0,2) -- node[pos=0,left] {$B$}(4,2);
\draw (4,0) -- node[pos=1,left] {$A$}(0,0);
\draw
(2.55,2.4)  
-- node[pos=0,right] {$C_{n-1}$}(1.9,1.5);
\draw
(1.7,1.8)    
--node[pos=1,right] {$C_{n-2}$} (2.6,1.4);
\node (dots) at (2.6,1.15) {$\vdots$};
\draw
(1.8,0.75)
--node[pos=1,right] {$C_2$}(2.60,0.3);
\draw
(2.4,0.6)
-- node[pos=1,right] {$C_1$}(2.1,-0.4);   
\end{tikzpicture}
\end{minipage}
\hfill
\begin{minipage}[c]{0.45\textwidth}
\centering
\begin{tikzpicture}[line width=1pt]
\draw (0,2)--(6,2);
\fill (0,2) circle (2pt);
\node[above] at (0,2) {$A$};
\fill (1.2,2) circle (2pt);
\node[above] at (1.2,2) {$C_1$};
\fill (2,2.2) circle (1.5pt);
\fill (2.5,2.2) circle (1.5pt);
\fill (3,2.2) circle (1.5pt);
\fill (3.5,2.2) circle (1.5pt);
\fill (4.8,2) circle (2pt);
\node[above] at (4.8,2) {$C_{n-1}$};
\fill (6,2) circle (2pt);
\node[above] at (6,2) {$B$};
\end{tikzpicture}   
\end{minipage}
\caption{The special fiber of reduction type $\mathrm{II}(n)$ and its dual graph}
    \label{fig:II}
\end{figure}

Let $C$ be a curve of genus $2$ with reduction type $\mathrm{II}(n)$. By the definition of reduction type $\mathrm{II}(n)$, the components of the special fiber of its stable model are two curves of genus $1$ meeting at a node of thickness $n$. The special fiber of the minimal regular model consists of two components of genus $1$ and a chain of $n-1$ rational components after the desingularization of the nodal point, see Figure\ref{fig:II}.

The following lemma gives an explicit formula for the Green function on vertices. 
\begin{lemma}\cite[Equation 4]{KayaComputation}\label{kayagreen}
Let $\Gamma$ be a metrized graph, and let $D=\sum_i a_i\Gamma_i$ be a divisor on the graph of $\deg(D)\ne-2$. For any two vertices $\Gamma_1$ and $\Gamma_2$ in $V(\Gamma)$, we obtain the following formula for the Green function associated to $D$
     \[G_D(\Gamma_1,\Gamma_2)=\frac{1}{\deg(D)+2}\left(\sum_i a_i j_{\Gamma_i}(\Gamma_1,\Gamma_2)-r(\Gamma_1,\Gamma_2)\right)+c_{D}\] where $c_{D}$ is a constant associated to $D$. Note that $j_{z}(x,y)$ is the voltage function and $r(x,y)\coloneq j_{y}(x,x)$ is the resistance function defined in \cite[\S2.2]{KayaComputation}. They satisfy the formula $j_{\Gamma_i}(\Gamma_1,\Gamma_2)=\frac{1}{2}(r(\Gamma_i,\Gamma_1)+r(\Gamma_i,\Gamma_2)-r(\Gamma_1,\Gamma_2))$.   
\end{lemma}
For a curve $C$ with semistable reduction, the resistance function $r(\,,\,)$ on $G_C$ can be computed as follows, see \cite[Remark~8.4]{MullerCanonical}. Let $M$ be the intersection matrix of the irreducible components with respect to an ordering of $(\Gamma_1,\ldots, \Gamma_n)$, and let $L=-M$. Let $v$ be the vector with entries zero except that the entry of $\Gamma_i$ is $1$, and the entry of $\Gamma_j$ is $-1$. The linear equation $Mx=v$ has solutions since the rank of $M$ is $n-1$ and $v\in \operatorname{Im}(M)$. Any two solutions differ by a constant vector, hence $-x^\top v$ is independent of the choice of the solution. We have the following formula \[r(\Gamma_i,\Gamma_j)=-x^\top v.\]
The matrix $L$ is the matrix of the discrete Laplacian operator on the graph, see \cite[\S2]{CinkirExplicitGreenfun}. We can view the solution $x$ as a piecewise linear function along the edges of $\Gamma$ except for $\Gamma_i$ and $\Gamma_j$ because $x$ is a harmonic function away from $\Gamma_i$ and $\Gamma_j$ and vertices of degree greater than $2$ on the graph, see \cite[\S2]{Chinburgcapacity}. This property of $x$ will be used in the following computation of the resistance function. 

\subsubsection{Computations for hyperelliptic curves with reduction type $\mathrm{III}(n)$}
Let $C$ be a hyperelliptic curve of genus $g$ over $K$ with semistable reduction type $\mathrm{III}(n)$. The special fiber of the minimal regular model consists of an irreducible component of genus $g-1$ labeled $A$, and a chain of $n-1$ rational components which are labeled $B_i$ such that $B_1$ and $B_{n-1}$ meet $A$. The resistance function is computed as follows.

\begin{lemma}\label{IIIa}
Let $C/K$ be a hyperelliptic curve of genus $g$ with semistable reduction type $\mathrm{III}(n)$. We have 
\begin{equation*}
    \begin{aligned}
        &r(\Gamma_A,\Gamma_{B_i})=\frac{i(n-i)}{n}\\
        &r(\Gamma_{B_i},\Gamma_{B_j})=|i-j|-\frac{(i-j)^2}{n}.
    \end{aligned}
\end{equation*}
where $\Gamma_A$ and $\Gamma_{B_i}$ are the respective vertices corresponding to the irreducible components $A$ and $B_i$.
\end{lemma}
\begin{proof}
The intersection matrix $M$ is given by 
\[\begin{pmatrix}
-2&1&0&\cdots&0&1\\
1&-2&1&\cdots&0&0\\
\vdots&\vdots&\vdots&\vdots&\vdots&\vdots\\
1&0&0&\cdots&1&-2\\
\end{pmatrix}\]
using the order $(\Gamma_A,\Gamma_{B_1},\ldots,\Gamma_{B_{n-1}})$. In order to compute $r(\Gamma_A,\Gamma_{B_i})$, we have to solve the linear equation $Mx=v$, where $v_1=1$, $v_{i+1}=-1$ and $v_j=0$ for $j\ne1,i+1$. Since the solution $x$ can be viewed as a piecewise linear function along the edges of $\Gamma$ except for $\Gamma_i$ and $\Gamma_j$, we obtain the formulas $x_1=0$ and $x_j=(j-1)c$ for $j\le i+1$ and $x_{i+2}=(i+1)c+1$, $x_j=x_{i+2}+(j-i-2)b$ for $j>i+2$ where $b$ and $c$ are constants to be determined. Applying these identities in the equation $Mx=v$, we obtain the formula $c=\frac{i-n}{n}$ and $b=\frac{i}{n}$. Hence we have the formula $r(\Gamma_A,\Gamma_{B_i})=-x^\top v=\frac{i(n-i)}{n}$. Moreover, $r(\Gamma_{B_i},\Gamma_{B_j})$ can be computed by the same method.
\end{proof}
Alternatively, by viewing the dual graph $G_C$ as an electrical network, one can prove Lemma~\ref{IIIa} using the parallel resistance formula. The associated Green functions can be computed as follows.
\begin{lemma}\label{IIIagreen}
Let $C/K$ be a hyperelliptic curve of genus $g$ with semistable reduction type $\mathrm{III}(n)$. We have 
\begin{equation*}
\begin{aligned}
   &G(\Gamma_A,\Gamma_{B_i})=-\frac{i(n-i)}{2gn}+c\\
    &G(\Gamma_{B_i},\Gamma_{B_j})=\frac{g-1}{2g}\left(\frac{i(n-i)}{n}+\frac{j(n-j)}{n}\right)+\frac{1}{2}\left(\frac{(j-i)^2}{n}-(j-i)\right)+c,\text{ where }i\le j\\
    &G(\Gamma_{B_i},\Gamma_{B_i})=\frac{g-1}{g}\frac{i(n-i)}{n}+c\\
    &G(\Gamma_A,\Gamma_A)=c
\end{aligned}  
\end{equation*}
where $c$ is a constant independent of the vertices.
\end{lemma}
\begin{proof}
    Note that the canonical divisor $K_G$ is $(2g-2)\Gamma_A$. By applying Lemma~\ref{kayagreen} with $D=K_G$ and combining Lemma~\ref{kayagreen} with Lemma~\ref{IIIa}, we have the desired values of the Green functions. 
\end{proof}
\subsubsection{Computations for genus $2$ curves with reduction type $\mathrm{II}(n)$}
Let $C/K$ be a curve of genus $2$ with semistable reduction type $\mathrm{II}(n)$. The special fiber of the minimal regular model consists of two elliptic curves labeled $A,B$ which are connected by a chain of $n-1$ components of genus $0$ labeled $(C_1,\ldots, C_{n-1})$.

\begin{lemma}\label{IIa}
Let $C/K$ be a curve of genus $2$ with semistable reduction type $\mathrm{II}(n)$. We have 
\begin{equation*}
    \begin{aligned}
        &r(\Gamma_A,\Gamma_B)=n&&r(\Gamma_A,\Gamma_{C_i})=i\\
        &r(\Gamma_B,\Gamma_{C_i})=n-i&&r(\Gamma_{C_i},\Gamma_{C_j})=j-i,\text{ where } i\leq j.\\
    \end{aligned}
\end{equation*}
\end{lemma}
\begin{proof}
    The computation of the resistance function $r$ follows the same method as that used in the proof of Lemma~\ref{IIIa}. Moreover, we can also prove it using the parallel resistance formula.
\end{proof}
\begin{lemma}\label{IIagreen}
Let $C/K$ be a curve of genus $2$ with semistable reduction type $\mathrm{II}(n)$. We have 
\begin{equation*}
  \begin{aligned}
      &G(\Gamma_A,\Gamma_B)=-\frac{1}{4}n+c,&&G(\Gamma_{C_i},\Gamma_{C_j})=\frac{1}{4}n+\frac{1}{2}(i-j)+c,\quad i\leq j\\
      &G(\Gamma_B,\Gamma_{C_i})=\frac{1}{4}(2i-n)+c,&&G(\Gamma_A,\Gamma_{C_i})=\frac{1}{4}(n-2i)+c\\
      &G(\Gamma_A,\Gamma_A)=\frac{1}{4}n+c,&&G(\Gamma_B,\Gamma_B)=\frac{1}{4}n+c\\
      &G(\Gamma_{C_i},\Gamma_{C_i})=\frac{1}{4}n+c\\
  \end{aligned} 
\end{equation*}
where $c$ is a constant.
\end{lemma}
\begin{proof}
    The canonical divisor is $K=\Gamma_A+\Gamma_B$ by the definition in \S\ref{canonical}. Let $K$ be the divisor $D$ in Lemma~\ref{kayagreen}. The Green functions are computed by using Lemma~\ref{kayagreen} and Lemma~\ref{IIa}. 
\end{proof}
\begin{rem}
For any given configuration of the special fiber of a semistable regular model, the Green function can always be computed by the same method.
\end{rem}

\subsection{Explicit formulas for $T(a)$}\label{ch8.2}
\begin{defi}
Let $A$ be an abelian variety over a non-archimedean local field $K$ with valuation ring $R$. The \textit{N\'eron model} of $A$ is a smooth separated group scheme $\mathcal{A}$ over $\operatorname{Spec}R$ whose generic fiber is $A$, satisfying the following universal property, called the N\'eron mapping property:
 
 For every smooth $R$-scheme $\mathcal{X}$, a morphism $\mathcal{X}_K\to A$  extends uniquely to a morphism $\mathcal{X}\to\mathcal{A}$.
\end{defi}
\begin{ex}\cite[\S$\mathrm{IV}$~Theorem~6.1]{SilvermanAdvancedTopics}
    Let $E/K$ be an elliptic curve, and let $\mathcal{C}/R$ be a minimal proper regular model for $E/K$. The N\'eron model of $E/K$ is the smooth locus of $\mathcal{C}/R$.   
\end{ex}
\begin{defi}
Let $A/K$ be an abelian variety. Let $\mathcal{A}$ be the N\'eron model of $A/K$ and $\mathcal{A}_s$ be its special fiber. Denote by $\mathcal{A}^0\subseteq\mathcal{A}$ the open subgroup scheme whose fibers are the identity components, and let $\mathcal{A}^0_s=(\mathcal{A}^0)_s$ be the identity component of the special fiber. The \textit{group of components} of $\mathcal{A}$ is denoted by $\Phi=\mathcal{A}_s/\mathcal{A}_s^0$.      
\end{defi}

N\'eron \cite[\S11.5]{LangDiophantine} proved that for any divisor $D\in\operatorname{Div}(A)$, there is a map $\gamma_D\colon\Phi\to\mathbb{R}$ such that the N\'eron function associated to $D$ is expressed as \[\lambda_D(a)=i(a,D)+\gamma_D(\textbf{a})\] where $a\in A(K)\setminus|D|$ and $\textbf{a}$ is the specialization of $a$ in the component group $\Phi$. The intersection multiplicity $i(a,D)$ is defined in \cite[\S11.5]{LangDiophantine}.  

\begin{rem}\label{remark8.12}
Let $C$ be a curve of genus $2$. In \cite[Definition~3.1]{MullerCanonical}, a function $\mu$ on the Kummer surface of $J(C)$ is defined; it can be extended to a function on $J(C)$. It is used in the definition of canonical local heights on Kummer coordinates, see \cite[Definition~4.1]{MullerCanonical}. It is shown in \cite[\S8]{MullerCanonical} that N\'eron functions on $J(C)$ can be constructed using canonical local heights on Kummer coordinates. 

Assume that $C$ has semistable reduction and $\mu$ factors through the component group $\Phi$. This occurs, for example, when $C$ has semistable reduction type $\mathrm{III}(n)$. For $a=[P_1-P_2]\in J(C)(K)\setminus|\Theta|$, the equality $\mu(a)=-2T(a)+c$ is shown in the proof of \cite[Proposition~8.5]{MullerCanonical}, where $c$ is a constant independent of $a$. Furthermore, the equation $\mu(a)=r(\Gamma_{P_1},\Gamma_{P_2})$ is shown in \cite[Theorem~8.7]{MullerCanonical} where $\Gamma_{P_i}$ is the specialization of $P_i$ on $G(C)$. When $C$ has semistable reduction type $\mathrm{III}(n)$, a formula for $\mu$ is given in \cite[Proposition~9.1]{MullerCanonical}. When $C$ is a hyperelliptic curve of genus $3$ with semistable reduction type $\mathrm{III}(3)$, we will compute $T(a)$ and prove that $T(a)$ does not factor through $\Phi$.             
\end{rem}
\begin{ex}[N\'eron model of a Jacobian]\cite[Appendix~A]{MatthewSpecialization}
   Let $C/K$ be a curve with semistable reduction. Let $\mathcal{C}$ be a semistable regular model of $C$ over $R$. Let $I=\{\Gamma_1,...,\Gamma_n\}$ be the set of irreducible components of the geometric special fiber $\mathcal{C}_{\bar{k}}$. Let $J(C)$ be the Jacobian of $C$, and let $\mathcal{J}$ be the N\'eron model of $J(C)$ over $R$. The intersection pairing on $I=\{\Gamma_1,...,\Gamma_n\}$ induces a map as follows
  \begin{equation*}
      \begin{aligned}
          \alpha&\colon\mathbb{Z}^I\rightarrow\mathbb{Z}^I\\&(a_i)_i\mapsto(\sum_{\Gamma_j\in I}a_j(\Gamma_i\cdot \Gamma_j))_i
      \end{aligned}
  \end{equation*}
where $(\Gamma_i\cdot\Gamma_j)$ is the intersection multiplicity on the special fiber $\mathcal{C}_{\bar{k}}$. The degree map from $\mathbb{Z}^I$ to $\mathbb{Z}$ is defined as $\deg\colon(a_i)_{i\in I}\mapsto\sum a_i$. It is shown in \cite[\S9.6~Theorem~1]{Boschneronmodles} that there is a complex \[\mathbb{Z}^I\xrightarrow{\alpha}\mathbb{Z}^I\xrightarrow{\deg}\mathbb{Z}\] and the component group $
\Phi$ is canonically isomorphic to $\operatorname{Ker}(\deg)/\operatorname{Im}(\alpha)$. The interpretation of the isomorphism $f\colon\Phi\cong\operatorname{Ker}(\deg)/\operatorname{Im}(\alpha)$ is as follows. 

Let $a\in J(C)(K)$, and let $\textbf{a}$ be the specialization of $a$ in the component group. Let $D\in\operatorname{Div}^0(C)$ such that $a=[D]$. For a closed point $P$, let $\bar{P}$ be the Zariski closure of $P$ in $\mathcal{C}$. We define $\bar{D}$ by linearity for any divisor $D$. Then we have $f(\textbf{a})=[C_i\mapsto(C_i\cdot\bar{D})]$. Moreover, if $D$ is supported on $K^{\operatorname{ur}}$-rational points, the isomorphism $f$ can be described as follows. The divisor $D$ can be written as $D=\sum_{P_i\in C(K^{\operatorname{ur}})}n_iP_i$ with $\sum n_i=0$. The closure of each point $P\in C(K^{\operatorname{ur}})$ meets a unique irreducible component of $\mathcal{C}_k$, which we denote by $\Gamma_P$. Then we have $f(\operatorname{sp}([D]))=[\sum n_i\Gamma_{P_i}]$.       
\end{ex}
\begin{lemma}\label{IIIndual}
Let $C$ be a hyperelliptic curve of genus $g$ with semistable reduction type $\mathrm{III}(n)$, and let $P\in C(K^{\operatorname{ur}})$. If $P$ specializes to $\Gamma_A$, then $\iota P$ also specializes to $\Gamma_A$ in the dual graph. And if $P$ specializes to $\Gamma_{B_i}$, then $\iota P$ specializes to $\Gamma_{B_{n-i}}$ in the dual graph.    
\end{lemma}
\begin{proof}
A graph $G$ is called hyperelliptic if it admits an involution $\iota_G$ such that the quotient graph $G/\iota_G$ is a tree, see\cite[Definition~2.2.1]{AlexanderTamagawa}. It is shown in \cite[Lemma~2.2.4]{AlexanderTamagawa} that when $G$ is the dual graph of the geometric special fiber of a semistable hyperelliptic curve, then there is an involution $\iota_G$, induced by the extension of $\iota$ to the semistable minimal regular model $\mathcal{C}$, such that $(G,\iota_G)$ is hyperelliptic. Moreover, the specialization map $\operatorname{sp}
\colon C(K^{\operatorname{ur}})\to G$ satisfies $\operatorname{sp}\circ\iota=\iota_G\circ\operatorname{sp}$.

The automorphism induced by $\iota$ on the special fiber of the minimal regular model maps each irreducible component isomorphically onto another irreducible component. Hence it preserves the genus of each irreducible component and it must preserve $A$, which is the unique irreducible component of positive genus. Therefore, $\operatorname{sp}(P)=\operatorname{sp}(\iota P)$ whenever $P$ specializes to $\Gamma_A$. Suppose now that $\operatorname{sp}(P)=\Gamma_{B_i}$, where $\Gamma_{B_i}$ is the vertex corresponding to the component $B_i$. Since the vertex $\Gamma_A$ is fixed by $\iota_G$ and the quotient graph $G/\iota_G$ is a tree, the involution must act on the chain of rational components by reflection. In particular, $\iota_G(\Gamma_{B_i})=\Gamma_{B_{n-i}}$. Therefore, $\operatorname{sp}(\iota P)=\Gamma_{B_{n-i}}$ by the identity $\operatorname{sp}\circ\iota=\iota_G\circ\operatorname{sp}$.      
\end{proof} 

It is shown in \cite[\S9]{MullerCanonical} that when a genus $2$ curve $C$ has semistable reduction type $\mathrm{III}(n)$, $T$ factors through $\Phi$. It is natural to ask whether this still holds in higher genus. We show that this is not the case.

\begin{ex}\label{IIIaT}
Let $C$ be a hyperelliptic curve of genus $3$ with semistable reduction type $\mathrm{III}(3)$. The specialization of $\infty$ is $\Gamma_A$. Let $a=[Q_1+Q_2+Q_3-3\infty]\in J(C)(K)\setminus|\Theta|$ and $Q_i\in C(K^{\operatorname{ur}})$. Let $\Gamma_{Q_i}$ denote the vertex corresponding to the specialization of $Q_i$ in the dual graph. Then all possible values of $T(a)$ are given as follows. Here $G(\Gamma_A,\Gamma_A)=c$ is a constant. We have
\begin{equation*}
    \begin{aligned}
        &T(a)=6c,&&\text{when }\Gamma_{Q_1}=\Gamma_{B_1},\Gamma_{Q_2}=\Gamma_{B_2},\Gamma_{Q_3}=\Gamma_{A}\\
        &T(a)=\frac{2}{3}+6c,&&\text{when }\Gamma_{Q_1}=\Gamma_{B_1},\Gamma_{Q_2}=\Gamma_{B_1},\Gamma_{Q_3}=\Gamma_{B_2}\\
        &T(a)=-\frac{1}{3}+6c,&&\text{when }\Gamma_{Q_1}=\Gamma_{B_1},\Gamma_{Q_2}=\Gamma_{B_1},\Gamma_{Q_3}=\Gamma_{A}\\
        &T(a)=-\frac{1}{3}+6c,&&\text{when }\Gamma_{Q_1}=\Gamma_{B_2},\Gamma_{Q_2}=\Gamma_{B_2},\Gamma_{Q_3}=\Gamma_{A}\\
        &T(a)=\frac{2}{3}+6c,&&\text{when }\Gamma_{Q_1}=\Gamma_{B_2},\Gamma_{Q_2}=\Gamma_{B_1},\Gamma_{Q_3}=\Gamma_{B_2}\\
        &T(a)=-\frac{1}{3}+6c,&&\text{when }\Gamma_{Q_1}=\Gamma_{A},\Gamma_{Q_2}=\Gamma_{B_2},\Gamma_{Q_3}=\Gamma_{A}\\
        &T(a)=-\frac{1}{3}+6c,&&\text{when }\Gamma_{Q_1}=\Gamma_{A},\Gamma_{Q_2}=\Gamma_{B_1},\Gamma_{Q_3}=\Gamma_{A}\\
        &T(a)=6c,&&\text{when }\Gamma_{Q_1}=\Gamma_{Q_2}=\Gamma_{Q_3}\in\{\Gamma_{A},\Gamma_{B_1},\Gamma_{B_2}\}.\\
    \end{aligned}
\end{equation*}
It is easy to find these possible values of $T(a)$ using Lemma~\ref{IIIagreen} and Lemma~\ref{IIIndual}. The component group $\Phi\cong\mathbb{Z}/3\mathbb{Z}$ is generated by $[\Gamma_{B_1}-\Gamma_{A}]$, see \cite[Theorem~9.6.1]{Boschneronmodles}. Moreover, $[\Gamma_{B_i}-\Gamma_{A}]=i\cdot[\Gamma_{B_1}-\Gamma_{A}]$. In order to show that $T$ does not factor through $\Phi$, it suffices to find two points $a,b\in J(K)\setminus|\Theta|$ with the same image in $\Phi$ but $T(a)\ne T(b)$. 

Let $a=[Q_1+Q_2+Q_3-3\infty]$ such that $\Gamma_{Q_1}=\Gamma_{B_1},\Gamma_{Q_2}=\Gamma_{B_1},\Gamma_{Q_3}=\Gamma_{B_2}$. The image of $a$ is $[\Gamma_{B_1}-\Gamma_{A}]$ in $\Phi$ and $T(a)=\frac{2}{3}+6c$. Let $b=[Q^{'}_1+Q^{'}_2+Q^{'}_3-3\infty]$ such that $\Gamma_{Q^{'}_1}=\Gamma_{B_2},\Gamma_{Q^{'}_2}=\Gamma_{B_2},\Gamma_{Q^{'}_3}=\Gamma_{A}$. The image of $b$ is also $[\Gamma_{B_1}-\Gamma_{A}]$ in $\Phi$, but $T(b)=-\frac{1}{3}+6c\ne T(a)$. Hence, $T$ does not factor through $\Phi$.            
\end{ex}
\begin{ex}
Let $C$ be a curve of genus $2$ with semistable reduction type $\mathrm{III}(3)$. The specialization of $\infty$ is $\Gamma_A$. Let $a=[P_1-P_2]\in J(C)(K)\setminus|\Theta|$ such that $P_i\in C(K^{\operatorname{ur}})$. We obtain the following values of $T(a)=G(\Gamma_{P_1},\Gamma_{P_2})+G(\Gamma_{P_1}+\Gamma_{P_2},\Gamma_\infty)$
\begin{equation*}
    \begin{aligned}
        &T(a)=-\frac{1}{3}+3c,&&\text{when }\Gamma_{P_1}=\Gamma_{A},\Gamma_{P_2}=\Gamma_{B_1}\\
        &T(a)=-\frac{1}{3}+3c,&&\text{when }\Gamma_{P_1}=\Gamma_{A},\Gamma_{P_2}=\Gamma_{B_2}\\
        &T(a)=-\frac{1}{3}+3c,&&\text{when }\Gamma_{P_1}=\Gamma_{B_1},\Gamma_{P_2}=\Gamma_{B_2}\\
        &T(a)=3c,&&\text{when }\Gamma_{P_1}=\Gamma_{P_2}.\\
    \end{aligned}
\end{equation*}
We have the equation $\mu(a)=-2T(a)+d$ by Remark~\ref{remark8.12} where $d$ is a constant. 
\end{ex}
When $C$ is a curve of genus $2$ with semistable reduction type $\mathrm{II}(n)$, the component group $\Phi$ is trivial by \cite[Proposition~9.6.10]{Boschneronmodles}. We show that the function $T$ does not factor through $\Phi$.
\begin{ex}
Let $C$ be a curve of genus $2$ with semistable reduction type $\mathrm{II}(2)$. We assume that the specialization of the point $\infty$ is $\Gamma_A$. Let $a=[P_1-P_2]\in J(C)(K)\setminus|\Theta|$ and $P_i\in C(K^{\operatorname{ur}})$. Then $T(a)=G(\Gamma_{P_1},\Gamma_{P_2})+G(\Gamma_{P_1}+\Gamma_{P_2},\Gamma_\infty)$. For the following three configurations, we obtain   
\begin{equation*}
    \begin{aligned}
        &T(a)=\frac{1}{2}+3c,&&\text{when }\Gamma_{P_1}=\Gamma_{A},\Gamma_{P_2}=\Gamma_{C_1}\\
        &T(a)=-\frac{1}{2}+3c,&&\text{when }\Gamma_{P_1}=\Gamma_{B},\Gamma_{P_2}=\Gamma_{C_1}\\
        &T(a)=-\frac{1}{2}+3c,&&\text{when }\Gamma_{P_1}=\Gamma_{A},\Gamma_{P_2}=\Gamma_{B}\\
    \end{aligned}
\end{equation*}
where $c$ is the constant in Lemma~\ref{IIagreen}. The computation of $T(a)$ proceeds in the same way as in Example~\ref{IIIaT}. If $T$ factored through the component group, then $T$ would be constant, because the component group is trivial. This contradicts the different values computed above.  
\end{ex}
\begin{rem}
Our method for computing $T(a)$ applies to all semistable reduction types of arbitrary thickness. The analysis of all possible specializations of the points $P_i$ becomes increasingly complicated as the thickness grows. For simplicity, we present only the above examples. 
\end{rem}
\printbibliography 
\end{document}